\documentclass[11pt,fleqn,a4paper]{article}
\usepackage{amssymb,latexsym,amsmath,amsfonts}
\usepackage{graphicx}

    \advance\textheight by \topskip

    \numberwithin{equation}{section}

    \newcommand\PVint{\mathop{\setbox0\hbox{$\displaystyle\intop$}%
    \hskip0.2\wd0%
    \vcenter{\hrule width0.6\wd0height0.5pt depth0.5pt}%
    \hskip-0.8\wd0%
    }\mskip-\thinmuskip\intop\nolimits}

    \def\ds{\displaystyle}
    
    \def\Re{{\rm Re \,}}
    \def\Im{{\rm Im \,}}
    \def\norm1{\|_{\raisebox{-.5ex}{\scriptsize{\!1}}}}

    \newtheorem{theorem}{Theorem}[section]
    \newtheorem{lemma}[theorem]{Lemma}
    
    \newtheorem{proposition}[theorem]{Proposition}
    
    \newtheorem{Definition}[theorem]{Definition}
    
    \newtheorem{Remark}[theorem]{Remark}
    \newenvironment{remark}{\begin{Remark}\rm}{\end{Remark}}
    \newtheorem{Example}[theorem]{Example}

    \newenvironment{proof}%
    {\rm \trivlist \item[\hskip \labelsep{\bf Proof. }]}%
    {\hspace*{\fill}$\Box$\endtrivlist}
    \newenvironment{varproof}%
    {\rm \trivlist \item[\hskip \labelsep{\bf Proof}]}%
    {\hspace*{\fill}$\Box$\endtrivlist}

    \newenvironment{varconjecture}%
    {\rm \trivlist \item[\hskip \labelsep{\bf Conjecture}]}%
    \endtrivlist

\begin{document}

    \begin{center} \Large\bf
        Strong asymptotics of the recurrence coefficients of
        orthogonal polynomials
        associated to the generalized Jacobi weight
    \end{center}

    \

    \begin{center}
        \large
        M. Vanlessen\footnote{Research Assistant of the Fund for Scientific Research -- Flanders (Belgium)}\\
            \normalsize \em
            Department of Mathematics, Katholieke Universiteit Leuven, \\
            Celestijnenlaan 200 B, 3001 Leuven, Belgium \\
            \rm maarten.vanlessen@wis.kuleuven.ac.be
    \end{center}\ \\[1ex]

\begin{abstract}
    We study asymptotics of the recurrence coefficients of orthogonal polynomials
    associated to the generalized Jacobi
    weight, which is a weight function with a finite number of algebraic singularities
    on $[-1,1]$. The recurrence coefficients can be written in terms of
    the solution of the
    corresponding Riemann-Hilbert problem for orthogonal polynomials.
    Using the steepest descent method of Deift and Zhou, we analyze the Riemann-Hilbert problem, and
    obtain complete asymptotic
    expansions of the recurrence coefficients. We will determine
    explicitly the order $1/n$ terms in the expansions. A critical step in the analysis
    of the Riemann-Hilbert problem will be the local analysis around the
    algebraic singularities, for which we use
    Bessel functions of appropriate order.
\end{abstract}

\section{Introduction}

We consider the generalized Jacobi weight
\begin{equation}\label{generalized_Jacobi_weight}
    w(x)=(1-x)^\alpha(1+x)^\beta h(x)\prod_{\nu=1}^{n_0}|x-x_\nu|^{2\lambda_\nu},
    \qquad\mbox{for $x\in(-1,1)$},
\end{equation}
where $n_0$ is a fixed number, with
\[
    -1<x_1<x_2<\cdots<x_{n_0}<1,\qquad 2\lambda_\nu>-1,\lambda_\nu\neq 0,\qquad
    \alpha,\beta>-1,
\]
and with $h$ real analytic and strictly positive on $[-1,1]$. The
points $x_1,\ldots ,x_{n_0}$ are called the algebraic
singularities of the weight. Throughout the paper we use $x_0=-1$
and $x_{n_0+1}=1$, for notational convenience. All the moments of
$w$ exist so that we have a sequence of orthogonal polynomials.
Denote the $n$-th degree orthonormal polynomial with respect to
the generalized Jacobi weight by $p_n(z)=\gamma_n z^n+\cdots$,
where $\gamma_n>0$. These orthonormal polynomials satisfy a three
term recurrence relation
\[
    x p_n(x)=a_{n+1} p_{n+1}(x)+b_n p_n(x)+a_n p_{n-1}(x),
\]
and we will investigate the asymptotic behavior of the recurrence
coefficients $a_n$ and $b_n$ as $n\to\infty$. The generalized
Jacobi weight has been studied before from other points of view in
\cite{Badkov,ErdMagNev,Vertesi1,Vertesi2} among others.

\medskip

For the pure Jacobi weight $(1-x)^\alpha (1+x)^\beta$ exact
expressions are known for the associated recurrence coefficients
$a_n$ and $b_n$, see \cite{Chihara,Szego}. The asymptotic behavior
is given by
\[
    a_n=\frac{1}{2}+O(1/n^2),\qquad b_n=O(1/n^2),\qquad \mbox{as $n\to\infty$.}
\]
In a previous paper with Kuijlaars, Mclaughlin and Van Assche
\cite{KMVV}, we considered the modified Jacobi weight
$(1-x)^\alpha(1+x)^\beta h(x)$. There, we were able to obtain
complete asymptotic expansions of the associated recurrence
coefficients in powers of $1/n$. It turned out that, as for the
pure Jacobi weight, the order $1/n$ terms in the expansions
vanished. The asymptotic behavior of the recurrence coefficients
of orthogonal polynomials associated to the generalized Jacobi
weight (\ref{generalized_Jacobi_weight}) has been studied before
by Golinskii \cite{Golinskii}. He has proven that
\begin{equation}\label{resultGolinskii}
    a_n=\frac{1}{2}+O(1/n),\qquad
    b_n=O(1/n),\qquad\mbox{as $n\to\infty$.}
\end{equation}
In this paper we will give stronger asymptotics. We will prove
that the $O(1/n)$ terms in (\ref{resultGolinskii}) can be
developed into complete asymptotic expansions in powers of $1/n$.
Here, in contrast with the (modified) Jacobi weight, the order
$1/n$ terms in the expansions will not vanish and we will
determine an explicit expression for them.

\medskip

Our approach is based on the characterization of orthogonal
polynomials via a Riemann-Hilbert problem, due to Fokas, Its and
Kitaev \cite{FokasItsKitaev}, and on an application of the
steepest descent method for Riemann-Hilbert problems of Deift and
Zhou \cite{DeiftZhou}. We have already applied this technique to
the modified Jacobi weight \cite{KMVV}, and in our case, the
general scheme is the same. The main difference lies in the fact
that we now have to do a local analysis around the algebraic
singularities as well (not just only around the endpoints $\pm
1$), which will be done with (modified) Bessel functions of
appropriate order. In the present paper, we will emphasize the
construction of the local parametrix near the algebraic
singularities, which is new. It will turn out that the order $1/n$
terms in the expansions of the recurrence coefficients come from
this parametrix. The Riemann-Hilbert approach has been applied
before to orthogonal polynomials, see
\cite{BleherIts,DKMVZ1,DKMVZ2,DKMVZ3,KriecherbauerMcLaughlin,KM,KMVV}.
Our result is the following.

\begin{theorem}\label{Theorem: asymptotic expansion recurrence coefficients}
    The recurrence coefficients $a_n$ and $b_n$ of orthogonal polynomials
    associated to the generalized Jacobi weight
    {\rm (\ref{generalized_Jacobi_weight})} have a complete asymptotic expansion of the form
    \begin{equation}\label{Theorem: Expansion}
        a_n\sim\frac{1}{2}+\sum_{k=1}^\infty\frac{A_k(n)}{n^k},
        \qquad\qquad
        b_n\sim\sum_{k=1}^\infty\frac{B_k(n)}{n^k},
    \end{equation}
    as $n\to\infty$. The coefficients $A_k(n)$ and $B_k(n)$ are
    explicitly computable for every $k$, and the coefficients with the $1/n$ term
    in the expansions are given by
    \begin{equation}\label{Theorem: A1(n)}
        A_1(n)=-\frac{1}{2}\sum_{\nu=1}^{n_0}\lambda_\nu\sqrt{1-x_\nu^2}\cos\bigl(2n\arccos
        x_\nu-\Phi_\nu\bigr),
    \end{equation}
    \begin{equation}\label{Theorem: B1(n)}
        B_1(n)=-\sum_{\nu=1}^{n_0}\lambda_\nu\sqrt{1-x_\nu^2}\cos\bigl((2n+1)\arccos x_\nu-\Phi_\nu\bigr),
    \end{equation}
    where
    \begin{eqnarray}
        \nonumber
        \Phi_\nu &=&
            \left(\alpha+\lambda_\nu+\sum_{k=\nu+1}^{n_0}2\lambda_k\right)\pi
            -\left(\alpha+\beta+\sum_{k=1}^{n_0}2\lambda_k\right)\arccos
            x_\nu \\[2ex]
        \label{Theorem: Phi}
        &&  -\,
            \frac{\sqrt{1-x_\nu^2}}{\pi}\PVint_{-1}^1\frac{\log
            h(t)}{\sqrt{1-t^2}}\frac{dt}{t-x_\nu}.
    \end{eqnarray}
    The integral in {\rm (\ref{Theorem: Phi})} is a Cauchy principal value integral.
\end{theorem}

This theorem shows that $n(2a_n-1)$ and $nb_n$ are oscillatory and
asymptotically behave like a superposition of $n_0$ wave functions
$A_\nu\cos(\omega_\nu n-\phi_\nu)$ with amplitudes
$A_\nu=-\lambda_\nu\sqrt{1-x_\nu^2}$, frequencies
$\omega_\nu=2\arccos x_\nu$, and phase shifts $\phi_\nu$ which are
different for $n(2a_n-1)$ and $n b_n$. The amplitude $A_\nu$
depends on the location and the strength of the singularity
$x_\nu$, while the frequency $\omega_\nu$ depends only on the
location of $x_\nu$. The strengths of the other singularities has
influence on the phase shift $\phi_\nu$. This discussion shows
that the $O(1/n)$ behavior of the recurrence coefficients is
intimately related to the behavior of our weight near the
singularities. Note that if we have no singularities (i.e.\!
$\lambda_1=\cdots =\lambda_{n_0}=0$) all the amplitudes in the
wave functions vanish. This implies that the order $1/n$ terms in
the expansions of the recurrence coefficients vanish, which is in
agreement with the case of the modified Jacobi weight \cite{KMVV}.

\begin{remark}
    We have restricted ourselves to determine only the
    the order $1/n$ terms in the expansions of the recurrence
    coefficients. It is possible to determine the higher
    order terms in the same way if we work hard enough, but the
    calculations will be a mess.
\end{remark}

We will now compare our result with a conjecture of Magnus
\cite{Magnus} about the asymptotic behavior of the recurrence
coefficients of orthogonal polynomials associated to the weight
\begin{equation}\label{weight Magnus}
        w(x)=\left\{\begin{array}{ll}
            B(1-x)^\alpha (1+x)^\beta (x_1-x)^{2\lambda}, & \qquad
            \mbox{for $x\in(-1,x_1)$,}\\[1ex]
            A(1-x)^\alpha (1+x)^\beta (x-x_1)^{2\lambda}, & \qquad
            \mbox{for $x\in(x_1,1)$.}
        \end{array}\right.
\end{equation}
where $A$ and $B$ are positive constants, with $-1<x_1<1$, where
$\alpha,\beta>-1$ and where $2\lambda>-1$. This weight allows a
jump at $x_1$, and is of the form
(\ref{generalized_Jacobi_weight}) only if $A=B$. The conjecture is
the following.

\begin{varconjecture}\textbf{of Magnus \cite{Magnus}.}
    The recurrence coefficients of orthogonal polynomials
    associated to the weight (\ref{weight Magnus}) satisfy
    \begin{eqnarray}
        \label{Conjecture: an}
        a_n
        &=&
            \frac{1}{2}-\frac{M}{n}\cos\Bigl(2n\arccos x_1
            -4\mu\log (4n \sin\arccos
            x_1)-\Phi\Bigr)+o(1/n), \\[1ex]
        \label{Conjecture: bn}
        b_n
        &=&
            -\frac{2M}{n}\cos\Bigl((2n+1)\arccos x_1-4\mu\log (4n\sin\arccos
            x_1)-\Phi\Bigr)+o(1/n),
    \end{eqnarray}
    as $n\to\infty$. Here
    \begin{equation}\label{Conjecture: mu M}
        \mu=\frac{1}{2\pi}\log(B/A), \qquad M=\frac{1}{2}\sqrt{\lambda^2+\mu^2}\sqrt{1-x_1^2},
    \end{equation}
    \begin{equation}\label{Conjecture: Phi}
        \Phi=(\alpha+\lambda)\pi
            -(\alpha+\beta+2\lambda)\arccos x_1-2\arg
        \Gamma(\lambda+i\mu)-\arg(\lambda+i\mu).
    \end{equation}
\end{varconjecture}

\medskip

\noindent We want to show that, as a consequence of Theorem
\ref{Theorem: asymptotic expansion recurrence coefficients}, the
conjecture is true for the case $A=B$. To this end, we need to
reformulate the conjecture for this case. If $A=B$ we have by
(\ref{Conjecture: mu M}) and (\ref{Conjecture: Phi}) that
$\mu=0,\, M=(|\lambda|/2)\sqrt{1-x_1^2}$, and
\[
    \Phi=
    \left\{
    \begin{array}{ll}
        (\alpha+\lambda)\pi-(\alpha+\beta+2\lambda)\arccos
        x_1, \qquad \mbox{if $\lambda>0$,}\\[1ex]
        (\alpha+\lambda)\pi-(\alpha+\beta+2\lambda)\arccos
        x_1-3\pi, \qquad \mbox{if $\lambda<0$.}
    \end{array}\right.
\]
Inserting this into (\ref{Conjecture: an}) and (\ref{Conjecture:
bn}) the conjecture becomes.

\begin{varconjecture}\textbf{of Magnus for the case $A=B$.} The
recurrence coefficients of orthogonal polynomials associated to
the weight $A(1-x)^\alpha (1+x)^\beta |x-x_1|^{2\lambda}$ satisfy
\[
    a_n=\frac{1}{2}-\frac{\lambda}{2n}\sqrt{1-x_1^2}\cos\Bigl(2n\arccos
    x_1-\hat\Phi\Bigr)+o(1/n),
\]
\[
    b_n=-\frac{\lambda}{n}\sqrt{1-x_1^2}\cos\Bigl((2n+1)\arccos
    x_1-\hat\Phi\Bigr)+o(1/n),
\]
as $n\to\infty$, where
$\hat\Phi=(\alpha+\lambda)\pi-(\alpha+\beta+2\lambda)\arccos x_1$.
\end{varconjecture}

\medskip

\noindent If we apply Theorem \ref{Theorem: asymptotic expansion
recurrence coefficients} for the case $n_0=1$ and $h=A$, and using
the fact that, see \cite{Gakhov},
\[
    \PVint_{-1}^1\frac{1}{\sqrt{1-t^2}}\frac{dt}{t-x_1}=0,
\]
we see that the conjecture is true for the case $A=B$. We even
have more since we were able to obtain complete asymptotic
expansions of the recurrence coefficients, allow the analytic
factor $h$, and allow more singularities. The full conjecture
(i.e.\! $A\neq B$) remains open.

\medskip

We will also compare our result with Nevai's result \cite{Nevai}
for an even positive weight $\rho(x)|x|^{2\lambda}$ on $[-1,1]$
with $\rho$ and $\rho'$ continuous and with $2\lambda>-1$. Nevai
showed in \cite{Nevai} that the recurrence coefficient $a_n$ of
orthogonal polynomials associated to this weight satisfies,
\begin{equation}\label{NevaiTheorem}
    a_n=\frac{1}{2}+(-1)^{n+1}\frac{\lambda}{2n}+o(1/n),\qquad\mbox{as $n\to\infty$.}
\end{equation}
We apply Theorem \ref{Theorem: asymptotic expansion recurrence
coefficients} to the weight $w(x)=(1-x^2)^\alpha
h(x)|x|^{2\lambda_1}$, with $h$ even. So, $w$ is of the form of
Nevai's weight. Using the fact that $h$ is even and $x_1=0$, we
have
\[
    \PVint_{-1}^1\frac{\log h(t)}{\sqrt{1-t^2}}\frac{dt}{t-x_1}=0
\]
since the integrand is an odd function. Since $n_0=1$,
$\alpha=\beta$, and $\arccos x_1=\pi/2$ this implies by
(\ref{Theorem: Phi}) that the phase constant $\Phi_1$ vanishes.
From (\ref{Theorem: Expansion}) and (\ref{Theorem: A1(n)}) we then
have
\[
    a_n=\frac{1}{2}+(-1)^{n+1}\frac{\lambda_1}{2n}+O(1/n^2),\qquad
    \mbox{as $n\to\infty$.}
\]
This is in agreement with Nevai's result, see
(\ref{NevaiTheorem}). The error is stronger since we have
$O(1/n^2)$ instead of $o(1/n)$. However, here we are dealing with
an even weight of the form $\rho(x)|x|^{2\lambda}$ with $\rho$
analytic, and Nevai's weights also include cases where $\rho$ is
non-analytic. Nevai also showed that the error is $O(1/n^2)$ if
$\rho$ is constant. Note, since the phase constant $\Phi_1$
vanishes, that by (\ref{Theorem: B1(n)}) the order $1/n$ term in
the expansion of $b_n$ vanishes. This is in agreement with the
fact that $b_n=0$ for an even weight.

\medskip

The present paper is organized as follows. In Section 2 we
formulate the theory of orthogonal polynomials as a
Riemann-Hilbert (RH) problem  for $Y$. In Section 3 we do the
asymptotic analysis of this RH problem. There, we want to obtain,
via a series of transformations $Y\mapsto T\mapsto S\mapsto R$, a
RH problem for $R$ with a jump matrix close to the identity
matrix. Then, $R$ is also close to the identity matrix. For the
last transformation $S\mapsto R$ we have to do a local analysis
near the endpoints and near the algebraic singularities. The local
analysis near the endpoints has already been done in \cite{KMVV},
but near the algebraic singularities it is new and will be done in
Section 4. In the last section we determine a complete asymptotic
expansion of the jump matrix for $R$. As a result, we obtain a
complete asymptotic expansion of $R$, which will be used to prove
Theorem \ref{Theorem: asymptotic expansion recurrence
coefficients}.

\section{RH problem for $Y$}

In this section we will characterize the orthogonal polynomials
via a $2\times 2$ matrix valued RH problem. This characterization
is due to Fokas, Its and Kitaev \cite{FokasItsKitaev}. We will
also write down the recurrence coefficients $a_n$ and $b_n$ in
terms of the solution of this RH problem.

We seek a $2\times 2$ matrix valued function $Y(z)=Y(z;n,w)$ that
satisfies the following RH problem.

\subsubsection*{RH problem for \boldmath$Y$:}
\begin{enumerate}
    \item[(a)]
        $Y(z)$ is  analytic for $z\in\mathbb C \setminus [-1,1]$.
    \item[(b)]
        $Y$ possesses continuous boundary values for $x \in (-1,1)\setminus\{x_1,\ldots ,x_{n_0}\}$
        denoted by $Y_{+}(x)$ and $Y_{-}(x)$, where $Y_{+}(x)$ and $Y_{-}(x)$
        denote the limiting values of $Y(z')$ as $z'$ approaches $x$ from
        above and below, respectively, and
        \begin{equation}\label{RHPYb}
            Y_+(x) = Y_-(x)
            \begin{pmatrix}
                1 & w(x) \\
                0 & 1
            \end{pmatrix},
            \qquad\mbox{for $x \in (-1,1)\setminus\{x_1,\ldots ,x_{n_0}\}$.}
        \end{equation}
    \item[(c)]
        $Y(z)$ has the following asymptotic behavior at infinity:
        \begin{equation} \label{RHPYc}
            Y(z)= (I+ O (1/z))
            \begin{pmatrix}
                z^{n} & 0 \\
                0 & z^{-n}
            \end{pmatrix}, \qquad \mbox{as $z\to\infty$.}
        \end{equation}
    \item[(d)]
        $Y(z)$ has the following behavior near $z=1$:
        \begin{equation}\label{RHPYd}
            Y(z)=\left\{
            \begin{array}{cl}
                O\begin{pmatrix}
                    1 & |z-1|^{\alpha} \\
                    1 & |z-1|^{\alpha}
                \end{pmatrix},
                &\mbox{if $\alpha<0$,} \\[2ex]
                O\begin{pmatrix}
                    1 & \log|z-1| \\
                    1 & \log|z-1|
                \end{pmatrix},
                &\mbox{if $\alpha=0$,} \\[2ex]
                O\begin{pmatrix}
                    1 & 1 \\
                    1 & 1
                \end{pmatrix},
                &\mbox{if $\alpha>0$,}
            \end{array}\right.
        \end{equation}
        as $z \to 1$, $z \in \mathbb C \setminus [-1,1]$.
    \item[(e)]
        $Y(z)$ has the following behavior near  $z=-1$:
        \begin{equation} \label{RHPYe}
            Y(z)=\left\{
            \begin{array}{cl}
                O\begin{pmatrix}
                    1 & |z+1|^{\beta} \\
                    1 & |z+1|^{\beta}
                \end{pmatrix}, &\mbox{if $\beta<0$,} \\[2ex]
                O\begin{pmatrix}
                    1 & \log|z+1| \\
                    1 & \log|z+1|
                \end{pmatrix},
                &\mbox{if $\beta=0$,} \\[2ex]
                O\begin{pmatrix}
                    1 & 1 \\
                    1 & 1
                \end{pmatrix},
                &\mbox{if $\beta>0$,}
            \end{array}\right.
        \end{equation}
        as $z \to -1$, $z \in \mathbb C\setminus [-1,1]$.
    \item[(f)]
        $Y(z)$ has the following behavior near  $z=x_\nu$, for every $\nu=1,\ldots ,n_0$:
        \begin{equation} \label{RHPYf}
            Y(z)=\left\{
            \begin{array}{cl}
                O\begin{pmatrix}
                    1 & |z-x_\nu|^{2\lambda_\nu} \\
                    1 & |z-x_\nu|^{2\lambda_\nu}
                \end{pmatrix}, &\mbox{if $\lambda_\nu< 0$,} \\[2ex]
                O\begin{pmatrix}
                    1 & 1 \\
                    1 & 1
                \end{pmatrix},
                &\mbox{if $\lambda_\nu>0$,}
            \end{array}\right.
        \end{equation}
        as $z \to x_\nu$, $z \in \mathbb C\setminus [-1,1]$.
\end{enumerate}

\begin{remark}
    The $O$-terms in (\ref{RHPYd}), (\ref{RHPYe}) and (\ref{RHPYf}) are to be taken
    entrywise. So for example $Y(z)= O\begin{pmatrix} 1 & |z-1|^{\alpha} \\
    1 & |z-1|^{\alpha}  \end{pmatrix}$  means that
    $Y_{11}(z) = O(1)$, $Y_{12}(z) = O(|z-1|^{\alpha})$, etc.
\end{remark}

If we take care of the algebraic singularities $x_\nu$ of the
generalized Jacobi weight in the same way as of the endpoints $\pm
1$ in \cite[Section 2]{KMVV} we obtain the following theorem.

\begin{theorem}
    The RH problem for $Y$ has a unique solution $Y(z)=Y(z;n,w)$ given by,
    \begin{equation}\label{RHPYsolution}
        Y(z) =
        \begin{pmatrix}
            \pi_n(z) & \frac{1}{2\pi i} \int_{-1}^1  \frac{\pi_n(x) w(x)}{x-z}dx \\[2ex]
            -2\pi i \gamma_{n-1}^2 \pi_{n-1}(z) & -\gamma_{n-1}^2 \int_{-1}^1 \frac{\pi_{n-1}(x)w(x)}{x-z} dx
        \end{pmatrix},
    \end{equation}
    where $\pi_n$ is the monic polynomial of degree $n$ orthogonal
    with respect to the weight $w$ and with $\gamma_n$ the leading
    coefficient of the orthonormal polynomial $p_n$.
\end{theorem}

\medskip

The recurrence coefficients $a_n$ and $b_n$ can be written in
terms of $Y$, see \cite{Deift,DKMVZ2,FokasItsKitaev,KMVV}. It is
known \cite{Deift} that
\begin{equation}\label{an in Y}
    a_n^2=\lim_{z\to\infty}z^2 Y_{12}(z;n,w)Y_{21}(z;n,w),
\end{equation}
\begin{equation}\label{bn in Y}
    b_n=\lim_{z\to\infty}\left(z-Y_{11}(z;n+1,w)Y_{22}(z;n,w)\right).
\end{equation}
So, in order to determine the asymptotics of the recurrence
coefficients, we need to do an asymptotic analysis of the RH
problem for $Y$.

\section{Asymptotic analysis of the RH problem for $Y$}

In this section we will do the asymptotic analysis of the RH
problem for $Y$. The idea is to obtain, via a series of
transformations
\[
    Y\mapsto T\mapsto S\mapsto R,
\]
a RH problem for $R$ which is normalized at infinity (i.e.
$R(z)\to I$ as $z\to\infty$) and whose jump matrix is close to the
identity matrix. As a result, the solution of the RH problem for
$R$ is also close to the identity matrix, cf.\!
\cite{Deift,DKMVZ2}.

As mentioned in the introduction, we point out that the asymptotic
analysis is analogous as in the case of the modified Jacobi
weight, see \cite{KMVV}. The main differences, which come from the
algebraic singularities, are:
\begin{itemize}
    \item In every step we have to take care of
    the growth condition near the algebraic singularities, which was included in the RH problem
    for $Y$ to control the behavior near these points.
    \item In the second transformation $T\to S$ the lens will be
    opened going through the algebraic singularities.
    \item We have to do a local analysis around the algebraic
    singularities, not just only around the endpoints.
    This is a new and critical step in the analysis of the RH
    problem for $Y$, and is the most important difference with the case
    of the modified Jacobi weight. To emphasize this, the construction of the parametrix near the
    algebraic singularities will be done in a
    separate section.
\end{itemize}

\subsection{First transformation $Y\to T$}

We will first transform the RH problem for $Y$ into a RH problem
for $T$ whose solution is bounded at infinity, and whose jump
matrix has oscillatory diagonal entries. Let
$\varphi(z)=z+(z^2-1)^{1/2}$ be the conformal mapping that maps
$\mathbb{C}\setminus[-1,1]$ onto the exterior of the unit circle,
and define
\begin{equation}\label{T in function of Y}
    T(z)=2^{n\sigma_3}Y(z)\varphi(z)^{-n\sigma_3}, \qquad\mbox{for $z\in\mathbb{C}\setminus[-1,1]$,}
\end{equation}
where $\sigma_3= \left(\begin{smallmatrix}
    1 & 0 \\
    0 & -1
\end{smallmatrix}\right)$ is the Pauli matrix. Then, $T$ is the
unique solution of the following equivalent RH problem, cf.\!
\cite[Section 3]{KMVV}.

\subsubsection*{RH problem for \boldmath$T$:}

\begin{enumerate}
    \item[(a)]
        $T(z)$ is analytic for $z\in\mathbb{C}\setminus[-1,1]$.
    \item[(b)]
        $T(z)$ satisfies the following jump relation on $(-1,1)\setminus\{x_1,\ldots ,x_{n_0}\}$:
        \begin{equation}\label{RHPTb}
            T_{+}(x)=T_{-}(x)
            \begin{pmatrix}
                \varphi_{+}(x)^{-2n} & w(x) \\
                0 & \varphi_{-}(x)^{-2n}
            \end{pmatrix}, \qquad\mbox{for $x\in(-1,1)\setminus\{x_1,\ldots ,x_{n_0}\}$.}
        \end{equation}
    \item[(c)]
        $T(z)$ has the following behavior at infinity:
        \begin{equation}\label{RHPTc}
            T(z) = I + O(1/z),
            \qquad \mbox{as $z \to \infty$.}
        \end{equation}
    \item[(d)]
        $T(z)$ has the same behavior as $Y(z)$ as $z\to 1$, given
        by (\ref{RHPYd}).
    \item[(e)]
        $T(z)$ has the same behavior as $Y(z)$ as $z\to -1$, given
        by (\ref{RHPYe}).
    \item[(f)]
        $T(z)$
        has the same behavior as $Y(z)$ as $z\to x_\nu$, given by
        (\ref{RHPYf}), for every $\nu=1,\ldots ,n_0$.
\end{enumerate}

\begin{remark}
    Condition (c) states that the RH problem for $T$ is normalized at
    infinity. Since $|\varphi_\pm(x)|=1$ for $x\in (-1,1)$ we have by (\ref{RHPTb})
    oscillatory diagonal entries in the jump matrix for $T$.
\end{remark}

\subsection{Second transformation $T\to S$}
\label{Subsection: Secont transformation}

We use the steepest descent method for RH problems of Deift and
Zhou \cite{DeiftZhou} to remove the oscillatory behavior in
(\ref{RHPTb}). See \cite{Deift, DKMVZ3} for an introduction. The
idea is to deform the contour so that the oscillatory diagonal
entries in the jump matrix for $T$ are transformed into
exponentially decaying off-diagonal entries. We then arrive at an
equivalent RH problem for $S$ on a lens shaped contour, with jumps
that converge to the identity matrix on the lips of the lens, as
$n\to\infty$. This step is referred to as the \emph{opening of the
lens}.

Since $h$ is real analytic and strictly positive on $[-1,1]$,
there is a neighborhood $U$ of $[-1,1]$ so that $h$ has an
analytic continuation to $U$, and so that the real part of $h$ is
strictly positive on $U$. Hence, the factor
$(1-x)^\alpha(1+x)^{\beta}h(x)$ has a non-vanishing analytic
continuation to $z\in
U\setminus\left((-\infty,-1]\cup[1,\infty)\right)$, given by
\[
    (1-z)^\alpha(1+z)^\beta h(z),
\]
with principal branches of powers.

\begin{figure}[h]
    \center{\resizebox{6cm}{!}{\includegraphics{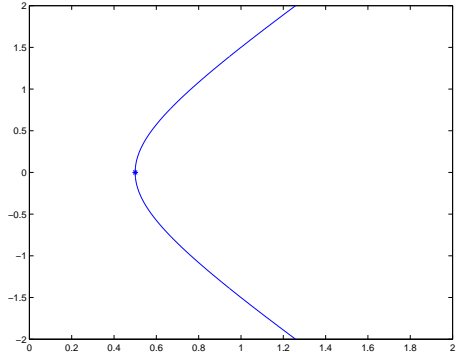}}
    \caption{The contour $\Gamma_{x_\nu}$ for the case $x_\nu=1/2$. The star marks $x_\nu$, and
    $K_{x_\nu}^l$ and $K_{x_\nu}^r$ are the sets of all points on
    the left, respectively right, of $\Gamma_{x_\nu}$.
    \label{figure1}}}
\end{figure}

To continuate the factor $|x-x_\nu|^{2\lambda_\nu}$ analytically,
where $\nu\in\{1,\ldots ,n_0\}$, we divide the complex plane into
two regions, which we denote by $K_{x_\nu}^l$ and $K_{x_\nu}^r$,
separated by a contour $\Gamma_{x_\nu}$ going through $x_\nu$, see
Figure \ref{figure1}. Here, $K_{x_\nu}^l$ and $K_{x_\nu}^r$ are
the sets of all points on the left, respectively right, of
$\Gamma_{x_\nu}$. We choose the contour $\Gamma_{x_\nu}$ so that
the images of $\Gamma_{x_\nu}\cap\mathbb{C}_+$ and
$\Gamma_{x_\nu}\cap\mathbb{C}_-$, under the mapping $\varphi$, are
the straight rays, restricted to the exterior of the unit circle,
with arguments $\arccos x_\nu$ and $-\arccos x_\nu$, respectively.
Here, $\mathbb{C}_+$ is used to denote the upper half-plane
$\{z\mid\Im z>0\}$, and $\mathbb{C}_-$ to denote the lower
half-plane $\{z\mid\Im z<0\}$. It turns out that $\Gamma_{x_\nu}$
is a hyperbola and goes vertically through $x_\nu$. We have made
an exact plot of $\Gamma_{x_\nu}$ for the case $x_\nu=1/2$, see
Figure \ref{figure1}. For $\nu=1,\ldots ,n_0$, the factor
$|x-x_\nu|^{2\lambda_\nu}$ has an analytic continuation to
$z\in\mathbb{C}\setminus\Gamma_{x_\nu}$, given by
\[
    \left\{\begin{array}{ll}
        (x_\nu-z)^{2\lambda_\nu},&\qquad \mbox{for $z\in K_{x_\nu}^l$,}\\[1ex]
        (z-x_\nu)^{2\lambda_\nu},&\qquad \mbox{for $z\in K_{x_\nu}^r$,}
    \end{array}\right.
\]
with again principal branches of powers.

\begin{remark}
    It seems a bit awkward to work with this choice of
    $\Gamma_{x_\nu}$ instead of with the vertical line going through $x_\nu$, but
    in Section \ref{Subsection: Construction Px01} this
    will become clear.
\end{remark}

As a result, the generalized Jacobi weight $w$, given by
(\ref{generalized_Jacobi_weight}), has a non-vanishing analytic
continuation to $z\in U\setminus\left((-\infty,
-1]\cup[1,\infty)\cup[\cup_{\nu=1}^{n_0}\Gamma_{x_\nu}]\right)$,
also denoted by $w$, given by
\begin{eqnarray}
    \nonumber
    w(z) &=&
        (1-z)^\alpha (1+z)^\beta h(z)\\[1ex]
    \label{Analytic extension of w}
    && \qquad\times\,
            \prod_{k=1}^{\nu}(z-x_k)^{2\lambda_k}\prod_{l=\nu+1}^{n_0}(x_l-z)^{2\lambda_l},
            \qquad\mbox{if $z\in K_{x_\nu}^r\cap K_{x_{\nu+1}}^l$,}
\end{eqnarray}
where $\nu=0,\ldots ,n_0$, and
$K_{x_0}^r=K_{x_{n_0+1}}^l=\mathbb{C}$.

\begin{remark}
    Note that only if $\lambda_\nu\in\mathbb{N}$ the analytic continuation of
    our weight is also analytic across the contour $\Gamma_{x_\nu}$.
\end{remark}

\begin{figure}[h]
    \center{\resizebox{12cm}{!}{\includegraphics{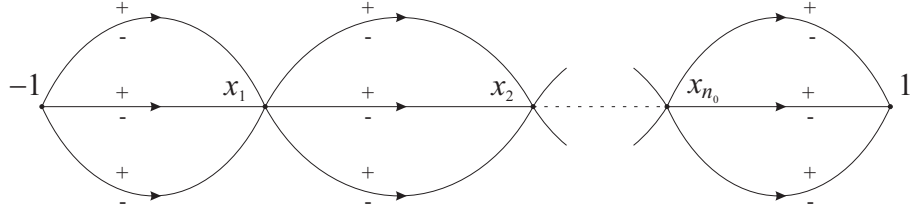}}
    \caption{The contour $\Sigma$.}\label{figure2}}
\end{figure}

The jump matrix (\ref{RHPTb}) for $T$ has the following
factorization into a product of three matrices, based on the fact
that $\varphi_+(x)\varphi_-(x)=1$ for $x\in(-1,1)$,
\begin{eqnarray}\nonumber
    \lefteqn{
    \begin{pmatrix}
        \varphi_{+}(x)^{-2n} & w(x) \\
        0 & \varphi_{-}(x)^{-2n}
    \end{pmatrix}} \\[1ex]
    & & =  \begin{pmatrix}
                1 & 0 \\
                w(x)^{-1}\varphi_{-}(x)^{-2n} & 1
            \end{pmatrix}
            \begin{pmatrix}
                0 & w(x) \\
                -w(x)^{-1} & 0
            \end{pmatrix}
            \begin{pmatrix}
                1 & 0 \\
                w(x)^{-1}\varphi_{+}(x)^{-2n} & 1
            \end{pmatrix}.
\end{eqnarray}
We note that $1/w$ does not have an analytic extension to a full
neighborhood of $(-1,1)$. Instead, it has an analytic continuation
to a neighborhood of $(x_\nu,x_{\nu+1})$, for every $\nu=0,\ldots
,n_0$, where $x_0=-1$ and $x_{n_0+1}=1$. We thus transform the RH
problem for $T$ into a RH problem for $S$ with jumps on the
oriented contour $\Sigma$, shown in Figure \ref{figure2}, that
goes through the algebraic singularities $x_\nu$. The precise form
of the lens $\Sigma$ will be determined in Section
\ref{Subsection: Construction Px01}. Of course it will be
contained in $U$. We write
\[
    \Sigma^o=\Sigma\setminus\{-1,x_1,\ldots ,x_{n_0},1\}.
\]

\noindent Let us define, as in \cite[Section 4]{KMVV},
\begin{equation} \label{S in function of T}
    S(z)=
    \left\{\begin{array}{cl}
        T(z), & \mbox{for $z$ outside the lens,} \\[2ex]
        T(z)
        \begin{pmatrix}
            1 & 0 \\
            -w(z)^{-1}\varphi(z)^{-2n} & 1
        \end{pmatrix}, & \mbox{for $z$ in the upper parts of the lens,} \\[2ex]
        T(z)
        \begin{pmatrix}
            1 & 0 \\
            w(z)^{-1}\varphi(z)^{-2n} & 1
        \end{pmatrix}, & \mbox{for $z$ in the lower parts of the lens.}
    \end{array} \right.
\end{equation}
Then, $S$ is the unique solution of the following equivalent RH
problem, cf.\! \cite[Section 4]{KMVV}.

\subsubsection*{RH problem for \boldmath$S$:}

\begin{enumerate}
    \item[(a)]
        $S(z)$ is analytic for $z\in\mathbb{C}\setminus\Sigma$.
    \item[(b)]
        $S(z)$ satisfies the following jump relations on $\Sigma^o$:
        \begin{equation}\label{RHPSb1}
            S_{+}(z)=S_{-}(z)
            \begin{pmatrix}
                1 & 0 \\
                w(z)^{-1}\varphi(z)^{-2n} & 1
            \end{pmatrix}, \qquad \mbox{for $z\in\Sigma^o\cap(\mathbb{C}_+\cup\mathbb{C}_-)$,}
        \end{equation}
        \begin{equation}\label{RHPSb2}
            S_{+}(x)=S_{-}(x)
            \begin{pmatrix}
                0 & w(x) \\
                -w(x)^{-1} & 0
            \end{pmatrix}, \qquad \mbox{for $x\in \Sigma^o\cap (-1,1)$.}
        \end{equation}
    \item[(c)]
        $S(z)$ has the following behavior at infinity:
        \begin{equation}\label{RHPSc}
            S(z)=I+O(1/z),\qquad\mbox{as $z\to\infty$.}
        \end{equation}
    \item[(d)]
        For $\alpha<0$, $S(z)$ has the following behavior as $z\to 1$:
        \begin{equation}\label{RHPSd1}
            S(z)=
            O\begin{pmatrix}
                1 & |z-1|^{\alpha} \\
                1 & |z-1|^{\alpha}
            \end{pmatrix}, \qquad\mbox{as $z\to 1, z\in\mathbb{C}\setminus\Sigma$.}
        \end{equation}
        For $\alpha=0$, $S(z)$ has the following behavior as $z\to 1$:
        \begin{equation}\label{RHPSd2}
            S(z)=
            O\begin{pmatrix}
                \log|z-1| & \log|z-1| \\
                \log|z-1| & \log|z-1|
            \end{pmatrix}, \qquad \mbox{as $z\to 1, z \in \mathbb C \setminus\Sigma$.}
        \end{equation}
        For $\alpha>0$, $S(z)$ has the following behavior as $z\to 1$:
        \begin{equation}\label{RHPSd3}
            S(z)=\left\{\begin{array}{cl}
                O\begin{pmatrix}
                    1 & 1 \\
                    1 & 1
                \end{pmatrix},& \mbox{as $z\rightarrow 1$ from outside the lens,} \\[2ex]
                O\begin{pmatrix}
                    |z-1|^{-\alpha} & 1 \\
                    |z-1|^{-\alpha} & 1
                \end{pmatrix}, & \mbox{as $z\to 1$ from inside the lens.}
            \end{array}\right.
        \end{equation}
    \item[(e)]
        $S(z)$ has the same behavior near $-1$ if we replace in (\ref{RHPSd1}),
        (\ref{RHPSd2}), and (\ref{RHPSd3}), $\alpha$ by $\beta$, $|z-1|$ by $|z+1|$ and take the
        limit $z \to -1$ instead of $z \to 1$.
    \item[(f)]
        For $\nu=1,\ldots ,n_0$, $S(z)$ has the following behavior
        as $z\to x_\nu$.
        For $\lambda_\nu<0$ we have
        \begin{equation}\label{RHPSf1}
            S(z)=
            O\begin{pmatrix}
                1 & |z-x_\nu|^{2\lambda_\nu} \\
                1 & |z-x_\nu|^{2\lambda_\nu}
            \end{pmatrix}, \qquad\mbox{as $z\to x_\nu, z\in\mathbb{C}\setminus\Sigma$.}
        \end{equation}
        For $\lambda_\nu>0$ we have
        \begin{equation}\label{RHPSf2}
            S(z)=\left\{\begin{array}{cl}
                O\begin{pmatrix}
                    1 & 1 \\
                    1 & 1
                \end{pmatrix},& \mbox{as $z\rightarrow x_\nu$ from outside the lens,} \\[2ex]
                O\begin{pmatrix}
                    |z-x_\nu|^{-2\lambda_\nu} & 1 \\
                    |z-x_\nu|^{-2\lambda_\nu} & 1
                \end{pmatrix}, & \mbox{as $z\to x_\nu$ from inside the lens.}
            \end{array}\right.
        \end{equation}
\end{enumerate}

Since $|\varphi(z)|>1$ for $z\in\mathbb{C}\setminus[-1,1]$ we see
from (\ref{RHPSb1}) that the oscillatory terms on the diagonal
entries in the jump matrix for $T$ have been transformed into
exponentially decaying off-diagonal entries in the jump matrix for
$S$ on the lips of the lens. So, the jump matrix for $S$ converges
exponentially fast to the identity matrix on the lips of the lens,
as $n\to\infty$. Hence, we expect that the leading order
asymptotics are determined by the solution of the following RH
problem.

\subsubsection*{RH problem for \boldmath$N$:}

\begin{enumerate}
    \item[(a)]
        $N(z)$ is analytic for $z \in \mathbb C \setminus [-1,1]$.
    \item[(b)]
        $N(z)$ satisfies the following jump relation on the interval $(-1,1)\setminus
        \{x_1,\ldots ,x_{n_0}\}$:
        \begin{equation} \label{RHPNb}
            N_+(x) = N_-(x)
            \begin{pmatrix}
                0 & w(x) \\
                -w(x)^{-1} & 0
            \end{pmatrix}, \qquad \mbox{for $x \in (-1,1)\setminus\{x_1,\ldots ,x_{n_0}\}$.}
        \end{equation}
    \item[(c)]
        $N(z)$ has the following behavior at infinity:
        \begin{equation} \label{RHPNc}
            N(z) = I + O(1/z), \qquad \mbox{as $z \to\infty$.}
        \end{equation}
\end{enumerate}
The solution of the RH problem for $N$ is referred to as the
parametrix for the outside region and it has been solved in
\cite[Section 5]{KMVV} using the Szeg\H{o} function associated
with the generalized Jacobi weight $w$,
\begin{equation}\label{Szegofunction}
    D(z)=\frac{(z-1)^{\alpha/2}(z+1)^{\beta/2}\prod_{\nu=1}^{n_0}(z-x_\nu)^{\lambda_\nu}}
        {\varphi(z)^{(\alpha+\beta+2\sum_{\nu=1}^{n_0}\lambda_\nu)/2}}\exp
        \left(\frac{(z^{2}-1)^{1/2}}{2\pi}\int_{-1}^{1}\frac{\log
        h(x)}{\sqrt{1-x^{2}}}\frac{dx}{z-x}\right).
\end{equation}
The Szeg\H{o} function $D(z)$ associated to $w$ is analytic and
non-zero for $z\in\mathbb{C}\setminus[-1,1]$, satisfies the jump
condition $D_+(x)D_-(x)=w(x)$ for $x\in(-1,1)\setminus\{x_1,\ldots
,x_{n_0}\}$, and $D_\infty=\lim_{z\to\infty}D(z)\in(0,+\infty)$.
The solution of the RH problem for $N$ is then given by, see
\cite[Section 5]{KMVV},
\begin{equation}\label{RHPNsolution}
    N(z) = D_{\infty}^{\sigma_3}
    \begin{pmatrix}
        \frac{a(z) + a(z)^{-1}}{2} & \frac{a(z) - a(z)^{-1}}{2i} \\[1ex]
        \frac{a(z) - a(z)^{-1}}{-2i} & \frac{a(z) + a(z)^{-1}}{2}
    \end{pmatrix} D(z)^{-\sigma_3},
\end{equation}
where
\begin{equation} \label{az}
    a(z) = \frac{(z-1)^{1/4}}{(z+1)^{1/4}}.
\end{equation}
For later use we have the following lemma.
\begin{lemma}\label{Lemma: D+}
    For every $\nu=0,\ldots ,n_0$,
    \begin{equation}
        D_+(x)=\sqrt{w(x)}e^{-i\psi_{\nu}(x)},\qquad\mbox{for $x_{\nu}<x<x_{\nu+1}$.}
    \end{equation}
    Here $x_0=-1$, $x_{n_0+1}=1$, and
    \begin{eqnarray}
        \nonumber
        \psi_\nu(x) &=&
            -\frac{1}{2}\left[\left(\alpha+\sum_{k=\nu+1}^{n_0}2\lambda_k\right)\pi
            -\left(\alpha+\beta+\sum_{k=1}^{n_0}2\lambda_k\right)\arccos x
            -\right. \\[1ex]
        \label{psi_nu}
        && \qquad\qquad \left.
            -\,
            \frac{\sqrt{1-x^2}}{\pi}
            \PVint_{-1}^1\frac{\log
            h(t)}{\sqrt{1-t^2}}\frac{dt}{t-x}\right].
    \end{eqnarray}
    The integral in {\rm (\ref{psi_nu})} is a Cauchy principal
    value integral.
\end{lemma}

\begin{proof}
    We rewrite the expression (\ref{Szegofunction}) for
    the Szeg\H{o} function as
    \begin{equation}\label{Szegofunctionlemma}
        D(z)=\frac{(z-1)^{\alpha/2}(z+1)^{\beta/2}\prod_{k=1}^{n_0}(z-x_k)^{\lambda_k}}
        {\varphi(z)^{(\alpha+\beta+2\sum_{k=1}^{n_0}\lambda_k)/2}}\exp
        \left(-i(z^2-1)^{1/2}\Phi(z)\right),
    \end{equation}
    where
    \[
        \Phi(z)=\frac{1}{2\pi i}\int_{-1}^{1}\frac{\log
        h(x)}{\sqrt{1-t^2}}\frac{dt}{t-z}.
    \]
    Now, we determine $D_+(x)$ for $x\in(x_\nu,x_{\nu+1})$. So, we need to take
    the $+$ boundary values for all quantities in (\ref{Szegofunctionlemma}). Using the Sokhotskii-Plemelj formula
    \cite[Section 4.2]{Gakhov}
    we have
    \[
        \Phi_+(x)=\frac{\log h(x)}{2\sqrt{1-x^2}}+\frac{1}{2\pi
        i}\PVint_{-1}^1\frac{\log
        h(t)}{\sqrt{1-t^2}}\frac{dt}{t-x},
    \]
    where the integral is a Cauchy principal value integral,
    so that, by (\ref{Szegofunctionlemma}) and the fact that $\varphi_+(x)=\exp(i\arccos x)$,
    the lemma is proved after an easy calculation.
\end{proof}

Before we can to do the third transformation we have to be
careful, since the jump matrices for $S$ and $N$ are not uniformly
close to each other near the endpoints $\pm 1$ and near the
algebraic singularities $x_\nu$. Therefore, a local analysis near
these points is necessary. Near the endpoints this has already
been done in \cite[Section 6]{KMVV}.

\medskip

We have constructed in \cite[Section 6]{KMVV} a parametrix $P_1$
in the disk $U_{\delta,1}$ with radius $\delta>0$, sufficiently
small, and center 1. This is a matrix valued function in
$U_{\delta,1}$, that has the same jumps as $S$ on $\Sigma$, that
matches with $N$ on the boundary $\partial U_{\delta,1}$ of
$U_{\delta,1}$,
\begin{equation}\label{Matching condition P}
    P_1(z)N^{-1}(z)=I+O(1/n), \qquad \mbox{as $n\to\infty$,
    uniformly for $z\in \partial U_{\delta,1}$,}
\end{equation}
and that has the same behavior as $S(z)$ near $z=1$. The
parametrix $P_1$ is given in \cite[Section 6]{KMVV}, and is
constructed out of Bessel function of order $\alpha$. We note that
the scalar function $W$ in \cite[(6.27)]{KMVV}, because of the
extra factor $\prod_{\nu=1}^{n_0}|x-x_\nu|^{2\lambda_\nu}$ in the
generalized Jacobi weight, should have an extra factor
$\prod_{\nu=1}^{n_0}(z-x_\nu)^{2\lambda_\nu}$.

Similarly we have constructed in \cite[Section 6]{KMVV} a
parametrix $P_{-1}$ in the disk $U_{\delta,-1}$ with radius
$\delta>0$ and center $-1$. This is a matrix valued function in
$U_{\delta,-1}$ that has the same jumps as $S$ on $\Sigma$, that
matches with $N$ on $\partial U_{\delta,-1}$
\begin{equation}\label{Matching condition Ptilde}
    P_{-1}(z)N^{-1}(z)=I+O(1/n), \qquad \mbox{as $n\to\infty$,
    uniformly for $z\in \partial U_{\delta,-1}$,}
\end{equation}
and that has the same behavior as $S(z)$ near $z=-1$. The
parametrix is given in \cite[Section 6]{KMVV}, and is constructed
out of Bessel functions of order $\beta$. We note that the scalar
function $\tilde W$ in \cite[(6.52)]{KMVV} should have an extra
factor $\prod_{\nu=1}^{n_0}(x_\nu-z)^{2\lambda_\nu}$.

We also have to construct a local parametrix $P_{x_\nu}$ near the
algebraic singularities $x_\nu$. Let $U_{\delta,x_\nu}$ be the
disk, with center $x_\nu$ and radius $\delta>0$ so that the
closures of the disks $U_{\delta,x_0},\ldots ,U_{\delta,
x_{n_0+1}}$ don't intersect and so that all the disks lie in $U$.
The construction of the parametrix $P_{x_\nu}$ will be done in
Section \ref{Subsection: Parametrix near the singularity}. For
now, let us assume that we have a $2\times 2$ matrix valued
function $P_{x_\nu}$ with the same jumps as $S$, that matches with
$N$ on $\partial U_{\delta,x_\nu}$,
\begin{equation}\label{Matching condition Pxnu}
    P_{x_\nu}(z)N^{-1}(z)=I+O(1/n), \qquad
    \mbox{as $n\to\infty$, uniformly for $z\in \partial U_{\delta,x_\nu}$,}
\end{equation}
and that has the same behavior as $S(z)$ near $z=x_\nu$.

\subsection{Third transformation $S\to R$}

Using the parametrix for the outside region, the parametrices near
the endpoints, and the parametrices near the algebraic
singularities we do the final transformation. Let us define the
matrix valued function $R$ as
\begin{equation}\label{R in function of S}
    R(z)=\left\{
    \begin{array}{ll}
        S(z)N^{-1}(z),& \qquad \mbox{for $z\in\mathbb{C}\setminus
            (\Sigma\cup [\cup_{\nu=0}^{n_0+1}U_{\delta,x_\nu}])$,} \\[1ex]
        S(z)P_{x_\nu}^{-1}(z), & \qquad \mbox{for $z\in U_{\delta,x_\nu}\setminus\Sigma$,
            and $\nu=0,\ldots ,n_0+1$.}
    \end{array}\right.
\end{equation}

\begin{remark}\label{Remark: invers exists}
    Note that the inverses of the parametrices exist. For $N,
    P_{-1}$ and $P_1$ this was already known, see \cite{KMVV}. In the next section we will show that
    $P_{x_\nu}$ is also invertible.
\end{remark}

\begin{figure}
    \center{\resizebox{15cm}{!}{\includegraphics{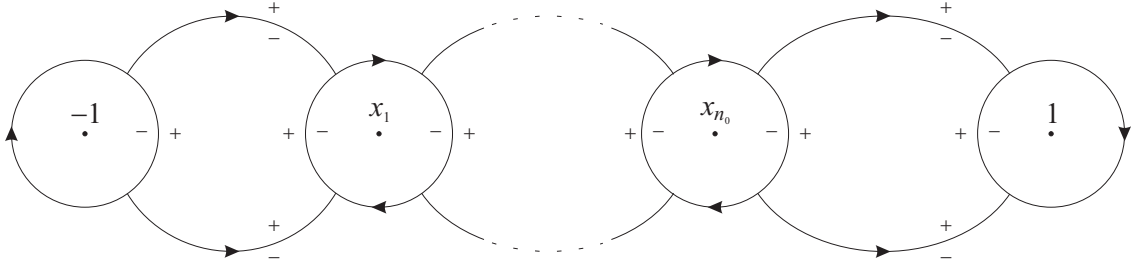}}
    \caption{The reduced system of contours $\Sigma_R$ with
    circles $U_{\delta,x_\nu}$ of radius $\delta$ and center $x_\nu$.}\label{figure3}}
\end{figure}

If we take care of the behavior near the algebraic singularities
in the same way as near the endpoints, it turns out that $R$
satisfies the following RH problem, cf.\! \cite[Section 7]{KMVV},
with jumps on the reduced system of contours $\Sigma_R$, see
Figure \ref{figure3}.

\subsubsection*{RH problem for \boldmath$R$:}

\begin{enumerate}
    \item[(a)]
        $R(z)$ is analytic for $z \in \mathbb C \setminus \Sigma_R$.
    \item[(b)]
        $R(z)$ satisfies the following jump relations on $\Sigma_R$:
        \begin{eqnarray}
            \nonumber
            R_+(z)
            &=&
                R_-(z)P_{x_\nu}(z) N^{-1}(z),  \\[2ex]
            \label{RHPRb1}
            & &
                \qquad \qquad \qquad
                \mbox{for $z\in \partial U_{\delta,x_\nu}$,
                    and $\nu=0,\ldots ,n_0+1$,} \\[2ex]
            \nonumber
            R_+(z)
            &=&
                R_-(z)N(z)
                \begin{pmatrix}
                    1 & 0 \\
                    w(z)^{-1}\varphi(z)^{-2n} & 1
                \end{pmatrix}N^{-1}(z), \\[2ex]
            \label{RHPRb2}
            & &
                \qquad \qquad \qquad
                \mbox{ for $z \in \Sigma_R \setminus(\cup_{\nu=0}^{n_0+1}\partial U_{\delta,x_\nu})$.}
        \end{eqnarray}
    \item[(c)]
        $R(z)$ has the following behavior at infinity:
        \begin{equation}\label{RHPRc}
            R(z)=I+O(1/z),\qquad\mbox{as $z\to\infty$.}
        \end{equation}
\end{enumerate}

By (\ref{Matching condition P}), (\ref{Matching condition Ptilde})
and (\ref{Matching condition Pxnu}), the jump matrices on the
circles are uniformly close to the identity matrix as
$n\to\infty$. On the lips of the lens we have by (\ref{RHPRb2}),
as in \cite[section 7]{KMVV}, that the jump matrix converges
uniformly to the identity matrix at an exponential rate. So, all
jump matrices are uniformly close to the identity matrix. This
implies that, cf.\! \cite{Deift,DKMVZ2},
\begin{equation}\label{R=I+O(1/n)}
    R(z)=I+O(1/n),\qquad \mbox{as $n\to\infty$, uniformly for $z\in\mathbb{C}\setminus\Sigma_R$.}
\end{equation}

\begin{remark}
    We will show in Section 5 that the $O(1/n)$ term in
    (\ref{R=I+O(1/n)}) can be developed into a complete asymptotic expansion in powers of
    $1/n$. This expansion will be used to prove Theorem \ref{Theorem: asymptotic expansion recurrence
    coefficients}.
\end{remark}

\section{Parametrix near the algebraic singularity $x_\nu$}
\label{Subsection: Parametrix near the singularity}

Fix $\nu\in\{1,\ldots ,n_0\}$. In this section we construct a
$2\times 2$ matrix valued function $P_{x_\nu}$ that satisfies the
following RH problem.

\subsubsection*{RH problem for \boldmath$P_{x_\nu}$:}
\begin{enumerate}
\item[(a)]
    $P_{x_\nu}(z)$ is defined and analytic for $z \in
    U_{\delta_0,x_\nu}\setminus\Sigma$ for some $\delta_0 > \delta$.
\item[(b)]
    $P_{x_\nu}(z)$ satisfies the following jump relations on
    $\Sigma^o\cap U_{\delta,x_\nu}$:
    \begin{equation}\label{RHPPx0b1}
        P_{x_\nu,+}(z) = P_{x_\nu,-}(z)
        \begin{pmatrix}
            1 & 0 \\
            w(z)^{-1}\varphi(z)^{-2n} & 1
        \end{pmatrix},
        \qquad \mbox{for $z\in \left(\Sigma^o \cap\mathbb{C}_\pm\right) \cap U_{\delta,x_\nu}$,}
    \end{equation}
    \begin{equation}\label{RHPPx0b2}
        P_{x_\nu,+}(x) = P_{x_\nu,-}(x)
        \begin{pmatrix}
            0 & w(x) \\
            -w(x)^{-1} & 0
        \end{pmatrix},
        \qquad \mbox{for $x \in \left(\Sigma^o\cap(-1,1)\right) \cap U_{\delta,x_\nu}$.}
    \end{equation}
\item[(c)]
    On $\partial U_{\delta,x_\nu}$ we have, as $n \to \infty$
    \begin{equation}\label{RHPPx0c}
        P_{x_\nu}(z) N^{-1}(z) = I + O(1/n),
        \qquad \mbox{uniformly for $z \in \partial U_{\delta,x_\nu}\setminus\Sigma$.}
    \end{equation}
\item[(d)]
    For $\lambda_\nu<0$, $P_{x_\nu}(z)$ has the following behavior as $z\to x_\nu$:
    \begin{equation}\label{RHPPx0d1}
        P_{x_\nu}(z)=
        O\begin{pmatrix}
            1 & |z-x_\nu|^{2\lambda_\nu} \\
            1 & |z-x_\nu|^{2\lambda_\nu}
        \end{pmatrix},
        \qquad \mbox{as $z\to x_\nu$.}
    \end{equation}
    For $\lambda_\nu>0$, $P_{x_\nu}(z)$ has the following behavior as $z\to x_\nu$:
    \begin{equation}\label{RHPPx0d2}
        P_{x_\nu}(z)=
        \left\{\begin{array}{cl}
            O\begin{pmatrix}
                1 & 1 \\
                1 & 1
            \end{pmatrix},
            & \mbox{as $z\to x_\nu$ from outside the lens,}\\[3ex]
            O\begin{pmatrix}
                |z-x_\nu|^{-2\lambda_\nu} & 1 \\
                |z-x_\nu|^{-2\lambda_\nu} & 1
            \end{pmatrix},
            & \mbox{as $z\to x_\nu$ from inside the lens.}
        \end{array}\right.
    \end{equation}
\end{enumerate}

We will work as follows. First, we construct a matrix valued
function that satisfies conditions (a), (b) and (d) of the RH
problem for $P_{x_\nu}$. For this purpose, we will transform (in
Section \ref{Subsection: Constant jump matrices}) this RH problem
into a RH problem for $P_{x_\nu}^{(1)}$ with constant jump
matrices and construct (in Section \ref{Subsection: Construction
Px01}) a solution of the RH problem for $P_{x_\nu}^{(1)}$.
Afterwards, we will also consider (in Section \ref{Subsection:
Matching condition}) the matching condition (c) of the RH problem
for $P_{x_\nu}$.

\subsection{Transformation to a RH problem with constant jump
matrices}
    \label{Subsection: Constant jump matrices}

Since $h$ is analytic in $U$ with positive real part, the scalar
function
\begin{eqnarray}
    \nonumber
    W_{x_\nu}(z)
    &=&
        (1-z)^{\alpha/2}(1+z)^{\beta/2}h^{1/2}(z)\prod_{k=1}^{\nu-1}(z-x_k)^{\lambda_k}
        \prod_{l=\nu+1}^{n_0}(x_l-z)^{\lambda_l} \\[1ex]
    \label{Wx0}
    &&
        \qquad\times\,
        \left\{\begin{array}{ll}
            (z-x_\nu)^{\lambda_\nu}, & \qquad \mbox{for $z\in (K_{x_\nu}^l\cap U)\setminus \mathbb{R}$,} \\[1ex]
            (x_\nu-z)^{\lambda_\nu}, & \qquad \mbox{for $z\in (K_{x_\nu}^r\cap U)\setminus \mathbb{R}$,}
        \end{array}\right.
\end{eqnarray}
is defined and analytic for $z\in U\setminus(\mathbb{R}\cup
\Gamma_{x_\nu})$. Here, we recall that $K_{x_\nu}^l$ and
$K_{x_\nu}^r$ are the sets of all points on the left, respectively
right, of $\Gamma_{x_\nu}$. We seek $P_{x_\nu}$ in the form
\begin{equation}\label{P in function of Px01}
    P_{x_\nu}(z)=E_{n,x_\nu}(z)P_{x_\nu}^{(1)}(z)W_{x_\nu}(z)^{-\sigma_3}\varphi(z)^{-n\sigma_3},
\end{equation}
where the matrix valued function $E_{n,x_\nu}$ is analytic in a
neighborhood of $U_{\delta,x_\nu}$, and $E_{n,x_\nu}$ will be
determined (in Section \ref{Subsection: Matching condition}) so
that the matching condition (c) of the RH problem for $P_{x_\nu}$
is satisfied.

\begin{figure}[h]
    \center{\resizebox{6cm}{!}{\includegraphics{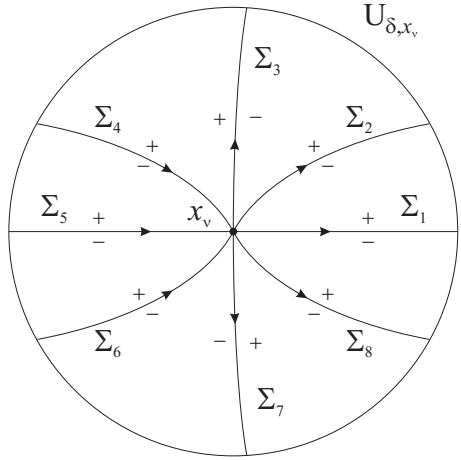}}
    \caption{The contour $\Sigma_{x_\nu}$. Here, $\Sigma_3\cup\Sigma_7$ is
    the part of $\Gamma_{x_\nu}$ in $U_{\delta,x_\nu}$, and the remainder is the part of
    $\Sigma$ in $U_{\delta,x_\nu}$.}\label{figure4}}
\end{figure}

Since $P_{x_\nu}$ has jumps on $\Sigma\cap U_{\delta,x_\nu}$, and
since $W_{x_\nu}$ has a jump on $\Gamma_{x_\nu}\cap U$, the matrix
valued function $P_{x_\nu}^{(1)}$ has jumps on the contour
\[
    \Sigma_{x_\nu}=(\Sigma\cup\Gamma_{x_\nu})\cap
    U_{\delta,x_\nu},
\]
see Figure \ref{figure4}. The contour $\Sigma_{x_\nu}$ consists of
8 parts, which we denote by $\Sigma_1,\ldots ,\Sigma_8$, as shown
in Figure \ref{figure4}. We write
\[
    \Sigma_{x_\nu}^o=\Sigma_{x_\nu}\setminus\{x_\nu\},\qquad
    \mbox{and}\qquad
    \Sigma_k^o=\Sigma_k\setminus\{x_\nu\},\quad\mbox{for $k=1,\ldots ,8$.}
\]

\medskip

In order to determine the jump matrices for $P_{x_\nu}^{(1)}$ we
need some information about the scalar function $W_{x_\nu}$. Write
\[
    \begin{array}{ll}
    K_{x_\nu}^{\rm I}=K_{x_\nu}^r \cap\mathbb{C}_+, & \qquad K_{x_\nu}^{\rm II}=K_{x_\nu}^l
    \cap\mathbb{C}_+, \\[1ex]
    K_{x_\nu}^{\rm III}=K_{x_\nu}^l \cap\mathbb{C}_-, & \qquad K_{x_\nu}^{\rm IV}=K_{x_\nu}^r
    \cap\mathbb{C}_-.
    \end{array}
\]
So, the sets $K_{x_\nu}^I,\ldots ,K_{x_\nu}^{IV}$ divide the
complex plane into four regions divided by the real axis and the
contour $\Gamma_{x_\nu}$. By (\ref{Analytic extension of w}) and
(\ref{Wx0}), we have for $z\in (K_{x_{\nu-1}}^r\cap
K_{x_{\nu+1}}^l)\cap U$,
\begin{equation}\label{Wx0^2}
    W_{x_\nu}^2(z)=
    \left\{\begin{array}{ll}
        w(z)e^{-2\pi i\lambda_\nu}, & \qquad \mbox{if
            $z\in K_{x_\nu}^{\rm I}\cup K_{x_\nu}^{\rm III}$,} \\[1ex]
        w(z)e^{2\pi i\lambda_\nu}, & \qquad \mbox{if
            $z\in K_{x_\nu}^{\rm II}\cup K_{x_\nu}^{\rm IV}$.}
    \end{array}\right.
\end{equation}
Here, we recall that $K_{x_0}^r=K_{x_{n_0+1}}^l=\mathbb{C}$. From
this we see that
\begin{equation}\label{Wx0+Wx0-}
    W_{x_\nu,+}(x)W_{x_\nu,-}(x)=w(x), \qquad \mbox{for $x\in(x_{\nu-1},x_{\nu+1})\setminus\{x_\nu\}$.}
\end{equation}
By (\ref{Wx0}) we have on the contour $\Gamma_{x_\nu}$,
\begin{equation}\label{Wx0+Wx0-invers}
    W_{x_\nu,+}(z)W_{x_\nu,-}(z)^{-1}=
        e^{\lambda_\nu\pi i}, \qquad\mbox{for $z\in\Gamma_{x_\nu}\cap U$.}
\end{equation}

We now have enough information about $W_{x_\nu}$ to determine the
jump matrices for $P_{x_\nu}^{(1)}$. First, we determine the jump
matrix on the lips of the lens. By (\ref{RHPPx0b1}), (\ref{P in
function of Px01}) and (\ref{Wx0^2}) the matrix valued function
$P_{x_\nu}^{(1)}$ should satisfy on
$\Sigma_2^o\cup\Sigma_4^o\cup\Sigma_6^o\cup\Sigma_8^o$ the jump
relation
\begin{eqnarray}
    \nonumber
    P_{x_\nu,+}^{(1)}(z)
    &=&
        P_{x_\nu,-}^{(1)}(z)\varphi(z)^{-n\sigma_3}W_{x_\nu}(z)^{-\sigma_3}
        \begin{pmatrix}
            1 & 0 \\
            w(z)^{-1}\varphi(z)^{-2n} & 1
        \end{pmatrix}
        W_{x_\nu}(z)^{\sigma_3}\varphi(z)^{n\sigma_3} \\[2ex]
    \nonumber
    &=&
        P_{x_\nu,-}^{(1)}(z)
        \begin{pmatrix}
            1 & 0 \\
            w(z)^{-1}W_{x_\nu}^2(z) & 1
        \end{pmatrix} \\[2ex]
    &=&
        P_{x_\nu,-}^{(1)}(z)
        \begin{pmatrix}
            1 & 0 \\
            e^{\pm 2\pi i\lambda_\nu} & 1
        \end{pmatrix},
\end{eqnarray}
where in $e^{\pm 2\pi i\lambda_\nu}$ the $+$ sign holds for
$z\in\Sigma_4^o\cup\Sigma_8^o$ and the $-$ sign for
$z\in\Sigma_2^o\cup\Sigma_6^o$.

Next, we determine the jump matrix on the interval. For $x\in
\Sigma_1^o\cup \Sigma_5^o$ we have by (\ref{RHPPx0b2}), (\ref{P in
function of Px01}), (\ref{Wx0+Wx0-}), and the fact that
$\varphi_+(x)\varphi_-(x)=1$,
\begin{eqnarray}
    \nonumber
    P_{x_\nu,+}^{(1)}(x)
    &=&
        P_{x_\nu,-}^{(1)}(x)W_{x_\nu,-}(x)^{-\sigma_3}\varphi_-(x)^{-n\sigma_3}
        \begin{pmatrix}
            0 & w(x) \\
            -w(x)^{-1} & 0
        \end{pmatrix}
        W_{x_\nu,+}(x)^{\sigma_3}\varphi_+(x)^{n\sigma_3} \\[2ex]
    \nonumber
    &=&
        P_{x_\nu,-}^{(1)}(x)
        \begin{pmatrix}
            0 & w(x)W_{x_\nu,+}(x)^{-1}W_{x_\nu,-}(x)^{-1} \\
            -w(x)^{-1}W_{x_\nu,+}(x)W_{x_\nu,-}(x) & 0
        \end{pmatrix} \\[2ex]
    &=&
        P_{x_\nu,-}^{(1)}(x)
        \begin{pmatrix}
            0 & 1 \\
            -1 & 0
        \end{pmatrix}.
\end{eqnarray}

And finally, we determine the jump matrix on the contour that goes
vertically through $x_\nu$. For $z\in \Sigma_3^o\cup\Sigma_7^o$ we
have by (\ref{P in function of Px01}) and (\ref{Wx0+Wx0-invers}),
\begin{equation}
    P_{x_\nu,+}^{(1)}(z)
    = P_{x_\nu,-}^{(1)}(z) W_{x_\nu,-}(z)^{-\sigma_3} W_{x_\nu,+}(z)^{\sigma_3}=
    P_{x_\nu,-}^{(1)}(z)e^{\lambda_\nu\pi i \sigma_3}.
\end{equation}

\medskip

We then see that we must look for a matrix valued function
$P_{x_\nu}^{(1)}$ that satisfies the following RH problem.

\subsubsection*{RH problem for \boldmath$P_{x_\nu}^{(1)}$:}
\begin{enumerate}
\item[(a)]
    $P_{x_\nu}^{(1)}(z)$ is defined and analytic for $z\in
    U_{\delta_0,x_\nu}\setminus(\Sigma\cup\Gamma_{x_\nu})$ for some $\delta_0 > \delta$.
\item[(b)]
    $P_{x_\nu}^{(1)}(z)$ satisfies the following jump relations on $\Sigma_{x_\nu}^o$:
    \begin{eqnarray}
        \label{RHPPx01b1}
        P_{x_\nu,+}^{(1)}(x)
        &=&
            P_{x_\nu,-}^{(1)}(x)
            \begin{pmatrix}
                0 & 1 \\
                -1 & 0
            \end{pmatrix},
            \qquad \mbox{for $x\in \Sigma_1^o\cup\Sigma_5^o$,} \\[2ex]
        \label{RHPPx01b2}
        P_{x_\nu,+}^{(1)}(z)
        &=&
            P_{x_\nu,-}^{(1)}(z)
            \begin{pmatrix}
                1 & 0 \\
                e^{-2\pi i\lambda_\nu} & 1
            \end{pmatrix},
            \qquad \mbox{for $z\in \Sigma_2^o\cup\Sigma_6^o$,} \\[2ex]
        \label{RHPPx01b3}
        P_{x_\nu,+}^{(1)}(z)
        &=&
            P_{x_\nu,-}^{(1)}(z) e^{\lambda_\nu\pi i\sigma_3},
            \qquad \mbox{for $z\in \Sigma_3^o\cup\Sigma_7^o$,} \\[2ex]
        \label{RHPPx01b4}
        P_{x_\nu,+}^{(1)}(z)
        &=&
            P_{x_\nu,-}^{(1)}(z)
            \begin{pmatrix}
                1 & 0 \\
                e^{2\pi i\lambda_\nu} & 1
            \end{pmatrix},
            \qquad \mbox{for $z\in \Sigma_4^o \cup\Sigma_8^o$.}
    \end{eqnarray}
\item[(c)]
    For $\lambda_\nu<0$, $P_{x_\nu}^{(1)}(z)$ has the following behavior as $z\to x_\nu$:
    \begin{equation}\label{RHPPx01c1}
        P_{x_\nu}^{(1)}(z)=
        O\begin{pmatrix}
            |z-x_\nu|^{\lambda_\nu} & |z-x_\nu|^{\lambda_\nu} \\
            |z-x_\nu|^{\lambda_\nu} & |z-x_\nu|^{\lambda_\nu}
        \end{pmatrix},
        \qquad \mbox{as $z\to x_\nu$.}
    \end{equation}
    For $\lambda_\nu>0$, $P_{x_\nu}^{(1)}(z)$ has the following behavior as $z\to x_\nu$:
    \begin{equation}\label{RHPPx01c2}
        P_{x_\nu}^{(1)}(z)=
        \left\{\begin{array}{cl}
            O\begin{pmatrix}
                |z-x_\nu|^{\lambda_\nu} & |z-x_\nu|^{-\lambda_\nu} \\
                |z-x_\nu|^{\lambda_\nu} & |z-x_\nu|^{-\lambda_\nu}
            \end{pmatrix},
            & \mbox{as $z\to x_\nu$ from outside the lens,} \\[3ex]
            O\begin{pmatrix}
                |z-x_\nu|^{-\lambda_\nu} & |z-x_\nu|^{-\lambda_\nu} \\
                |z-x_\nu|^{-\lambda_\nu} & |z-x_\nu|^{-\lambda_\nu}
            \end{pmatrix},
            & \mbox{as $z\to x_\nu$ from inside the lens.}
        \end{array}\right.
    \end{equation}
\end{enumerate}

\begin{remark}
    Condition (c) follows from condition (d) of
    the RH problem for $P_{x_\nu}$, since
    \[
        P_{x_\nu}(z)=E_{n,x_\nu}(z)P_{x_\nu}^{(1)}(z)W_{x_\nu}(z)^{-\sigma_3}\varphi(z)^{-n\sigma_3},
    \]
    where $\varphi(z)$ is bounded and bounded away from 0 near $z=x_\nu$, and where $W_{x_\nu}(z)$ behaves like
    $c|z-x_\nu|^{\lambda_\nu}$ as $z\to x_\nu$, with a non-zero constant $c$.
\end{remark}

\subsection{Construction of $P_{x_\nu}^{(1)}$}
    \label{Subsection: Construction Px01}

\begin{figure}
    \center{\resizebox{6cm}{!}{\includegraphics{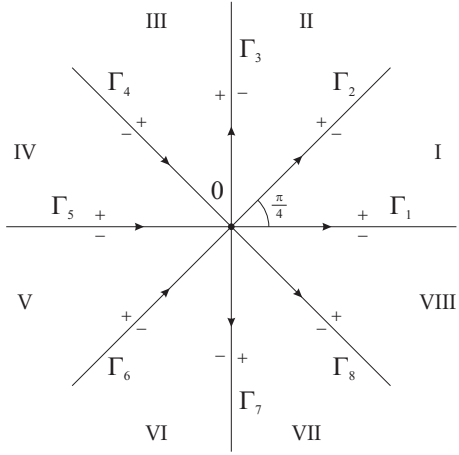}}
    \caption{The contour $\Sigma_\Psi$.}\label{figure5}}
\end{figure}

The construction of $P_{x_\nu}^{(1)}$ is based upon an auxiliary
RH problem for $\Psi_{\lambda}$ in the $\zeta$-plane with jumps on
the contour $\Sigma_\Psi=\cup_{i=1}^8 \Gamma_i$ consisting of
eight straight rays, oriented as in Figure \ref{figure5}, which
divides the complex plane into eight regions I--VIII, also shown
in Figure \ref{figure5}. We let $\lambda>-1/2$.

\subsubsection*{RH problem for \boldmath$\Psi_{\lambda}$:}
\begin{enumerate}
\item[(a)]
    $\Psi_{\lambda}(\zeta)$ is analytic for $\zeta\in\mathbb{C}\setminus\Sigma_\Psi$.
\item[(b)]
    $\Psi_{\lambda}(\zeta)$ satisfies the following jump relations on $\Sigma_\Psi$:
    \begin{eqnarray}
        \label{RHPPSIx0b1}
        \Psi_{\lambda,+}(\zeta)
        &=&
            \Psi_{\lambda,-}(\zeta)
            \begin{pmatrix}
                0 & 1 \\
                -1 & 0
            \end{pmatrix},
            \qquad \mbox{for $\zeta \in \Gamma_1\cup\Gamma_5$,} \\[2ex]
        \label{RHPPSIx0b2}
        \Psi_{\lambda,+}(\zeta)
        &=&
            \Psi_{\lambda,-}(\zeta)
            \begin{pmatrix}
                1 & 0 \\
                e^{-2\pi i \lambda}  & 1
            \end{pmatrix},
            \qquad \mbox{for $\zeta \in \Gamma_2\cup\Gamma_6$,} \\[2ex]
        \label{RHPPSIx0b3}
        \Psi_{\lambda,+}(\zeta)
        &=&
            \Psi_{\lambda,-}(\zeta) e^{\lambda \pi i \sigma_3},
            \qquad \mbox{for $\zeta \in \Gamma_3\cup\Gamma_7$,} \\[2ex]
        \label{RHPPSIx0b4}
        \Psi_{\lambda,+}(\zeta)
        &=&
            \Psi_{\lambda,-}(\zeta)
            \begin{pmatrix}
                1 & 0 \\
                e^{2\pi i \lambda}  & 1
            \end{pmatrix},
            \qquad \mbox{for $\zeta \in \Gamma_4\cup\Gamma_8$.}
    \end{eqnarray}
\item[(c)]
    For $\lambda<0$, $\Psi_{\lambda}(\zeta)$ has the following behavior as $\zeta\to 0$:
    \begin{equation}\label{RHPPSIx0c1}
        \Psi_{\lambda}(\zeta)=
        O\begin{pmatrix}
            |\zeta|^{\lambda} & |\zeta|^{\lambda} \\
            |\zeta|^{\lambda} & |\zeta|^{\lambda}
        \end{pmatrix},
        \qquad \mbox{as $\zeta\to 0$.}
    \end{equation}
    For $\lambda>0$, $\Psi_{\lambda}(\zeta)$ has the following behavior as $\zeta\to 0$:
    \begin{equation}\label{RHPPSIx0c2}
        \Psi_{\lambda_\nu}(\zeta)=
        \left\{\begin{array}{cl}
            O\begin{pmatrix}
                |\zeta|^{\lambda} & |\zeta|^{-\lambda} \\
                |\zeta|^{\lambda} & |\zeta|^{-\lambda}
            \end{pmatrix},
            & \mbox{as $\zeta\to 0$ for $\zeta\in$ II, III, VI, VII,} \\[3ex]
            O\begin{pmatrix}
                |\zeta|^{-\lambda} & |\zeta|^{-\lambda}\\
                |\zeta|^{-\lambda} & |\zeta|^{-\lambda}
            \end{pmatrix},
            & \mbox{as $\zeta\to 0$ for $\zeta\in$ I, IV, V, VIII.}
        \end{array}\right.
    \end{equation}
\end{enumerate}

We construct a solution $\Psi_{\lambda}$ of this RH problem out of
the modified Bessel functions $I_{\lambda\pm 1/2}$ and
$K_{\lambda\pm 1/2}$, and out of the Hankel functions
$H^{(1)}_{\lambda\pm 1/2}$ and $H^{(2)}_{\lambda\pm 1/2}$. For
$\zeta\in$ I, we define $\Psi_{\lambda}(\zeta)$ by
\begin{equation}\label{RHPPSIx0solution1}
    \Psi_{\lambda}(\zeta)= \frac{1}{2}\sqrt\pi\zeta^{1/2}
    \begin{pmatrix}
        H_{\lambda+\frac{1}{2}}^{(2)}(\zeta) &
            -i H_{\lambda+\frac{1}{2}}^{(1)}(\zeta) \\[2ex]
        H_{\lambda-\frac{1}{2}}^{(2)}(\zeta) &
            -i H_{\lambda-\frac{1}{2}}^{(1)}(\zeta)
    \end{pmatrix}
    e^{-(\lambda+\frac{1}{4}) \pi i \sigma_3}.
\end{equation}
For $\zeta\in$ II, by
\begin{equation}\label{RHPPSIx0solution2}
    \Psi_{\lambda}(\zeta)=
    \begin{pmatrix}
        \sqrt\pi\zeta^{1/2}I_{\lambda+\frac{1}{2}}(\zeta e^{-\frac{\pi i}{2}}) &
            -\frac{1}{\sqrt\pi}\zeta^{1/2}K_{\lambda+\frac{1}{2}}(\zeta e^{-\frac{\pi i}{2}}) \\[2ex]
        -i\sqrt\pi\zeta^{1/2}I_{\lambda-\frac{1}{2}}(\zeta e^{-\frac{\pi i}{2}}) &
            -\frac{i}{\sqrt\pi}\zeta^{1/2}K_{\lambda-\frac{1}{2}}(\zeta e^{-\frac{\pi i}{2}})
    \end{pmatrix}
    e^{-\frac{1}{2}\lambda\pi i\sigma_3}.
\end{equation}
For $\zeta\in$ III, by
\begin{equation}\label{RHPPSIx0solution3}
    \Psi_{\lambda}(\zeta)=
    \begin{pmatrix}
        \sqrt\pi\zeta^{1/2}I_{\lambda+\frac{1}{2}}(\zeta e^{-\frac{\pi i}{2}}) &
            -\frac{1}{\sqrt\pi}\zeta^{1/2}K_{\lambda+\frac{1}{2}}(\zeta e^{-\frac{\pi i}{2}}) \\[2ex]
        -i\sqrt\pi\zeta^{1/2}I_{\lambda-\frac{1}{2}}(\zeta e^{-\frac{\pi i}{2}}) &
            -\frac{i}{\sqrt\pi}\zeta^{1/2}K_{\lambda-\frac{1}{2}}(\zeta e^{-\frac{\pi i}{2}})
    \end{pmatrix}
    e^{\frac{1}{2}\lambda\pi i\sigma_3}.
\end{equation}
For $\zeta\in$ IV, by
\begin{equation}\label{RHPPSIx0solution4}
    \Psi_{\lambda}(\zeta)=\frac{1}{2}\sqrt\pi(-\zeta)^{1/2}
    \begin{pmatrix}
        i H_{\lambda+\frac{1}{2}}^{(1)}(-\zeta) &
            - H_{\lambda+\frac{1}{2}}^{(2)}(-\zeta) \\[2ex]
        -i H_{\lambda-\frac{1}{2}}^{(1)}(-\zeta) &
            H_{\lambda-\frac{1}{2}}^{(2)}(-\zeta)
    \end{pmatrix}
    e^{(\lambda+\frac{1}{4})\pi i \sigma_3}.
\end{equation}
For $\zeta\in$ V, by
\begin{equation}\label{RHPPSIx0solution5}
    \Psi_{\lambda}(\zeta)=\frac{1}{2}\sqrt\pi(-\zeta)^{1/2}
    \begin{pmatrix}
        - H_{\lambda+\frac{1}{2}}^{(2)}(-\zeta) &
            -i H_{\lambda+\frac{1}{2}}^{(1)}(-\zeta) \\[2ex]
        H_{\lambda-\frac{1}{2}}^{(2)}(-\zeta) &
            i H_{\lambda-\frac{1}{2}}^{(1)}(-\zeta)
    \end{pmatrix}
    e^{-(\lambda+\frac{1}{4})\pi i \sigma_3}.
\end{equation}
For $\zeta\in$ VI, by
\begin{equation}\label{RHPPSIx0solution6}
    \Psi_{\lambda}(\zeta)=
    \begin{pmatrix}
        -i\sqrt{\pi}\zeta^{1/2}I_{\lambda+\frac{1}{2}}(\zeta e^{\frac{\pi i}{2}}) &
            -\frac{i}{\sqrt\pi}\zeta^{1/2}K_{\lambda+\frac{1}{2}}(\zeta e^{\frac{\pi i}{2}}) \\[2ex]
        \sqrt{\pi}\zeta^{1/2}I_{\lambda-\frac{1}{2}}(\zeta e^{\frac{\pi i}{2}}) &
            -\frac{1}{\sqrt\pi}\zeta^{1/2}K_{\lambda-\frac{1}{2}}(\zeta e^{\frac{\pi i}{2}})
    \end{pmatrix}
    e^{-\frac{1}{2}\lambda\pi i \sigma_3}.
\end{equation}
For $\zeta\in$ VII, by
\begin{equation}\label{RHPPSIx0solution7}
    \Psi_{\lambda}(\zeta)=
    \begin{pmatrix}
        -i\sqrt{\pi}\zeta^{1/2}I_{\lambda+\frac{1}{2}}(\zeta e^{\frac{\pi i}{2}}) &
            -\frac{i}{\sqrt\pi}\zeta^{1/2}K_{\lambda+\frac{1}{2}}(\zeta e^{\frac{\pi i}{2}}) \\[2ex]
        \sqrt{\pi}\zeta^{1/2}I_{\lambda-\frac{1}{2}}(\zeta e^{\frac{\pi i}{2}}) &
            -\frac{1}{\sqrt\pi}\zeta^{1/2}K_{\lambda-\frac{1}{2}}(\zeta e^{\frac{\pi i}{2}})
    \end{pmatrix}
    e^{\frac{1}{2}\lambda\pi i \sigma_3}.
\end{equation}
And finally, for $\zeta\in$ VIII, we define it by
\begin{equation}\label{RHPPSIx0solution8}
    \Psi_{\lambda}(\zeta)= \frac{1}{2}\sqrt\pi\zeta^{1/2}
    \begin{pmatrix}
        -i H_{\lambda+\frac{1}{2}}^{(1)}(\zeta) &
            -H_{\lambda+\frac{1}{2}}^{(2)}(\zeta) \\[2ex]
        -iH_{\lambda-\frac{1}{2}}^{(1)}(\zeta) &
            -H_{\lambda-\frac{1}{2}}^{(2)}(\zeta)
    \end{pmatrix}
    e^{(\lambda+\frac{1}{4}) \pi i \sigma_3}.
\end{equation}

\begin{theorem}
    The matrix valued function $\Psi_{\lambda}$, defined by
    {\rm(\ref{RHPPSIx0solution1})--(\ref{RHPPSIx0solution8})},
    is a solution of the RH problem for $\Psi_{\lambda}$.
\end{theorem}
\begin{proof}
    The functions $I_{\lambda\pm 1/2},K_{\lambda\pm 1/2}, H_{\lambda\pm 1/2}^{(1)}$
    and $H_{\lambda\pm 1/2}^{(2)}$
    are defined and analytic in the complex plane with a
    branch cut along the negative real axis. So, the matrix valued
    function $\Psi_{\lambda}$ defined by (\ref{RHPPSIx0solution1})--(\ref{RHPPSIx0solution8})
    is analytic in the respective regions, and condition (a) of the RH problem is therefore satisfied.
    Condition (c) follows easily from
    \cite[formulas 9.1.9, 9.6.7 and 9.6.9]{AbramowitzStegun}. So, it remains to prove
    that jump conditions (\ref{RHPPSIx0b1})--(\ref{RHPPSIx0b4}) are satisfied.

    \medskip

    \textit{Jump conditions {\rm (\ref{RHPPSIx0b1})} and {\rm (\ref{RHPPSIx0b3})}}:
    By inspection, it is easy to see that these
    jump conditions are satisfied.

    \medskip

    \textit{Jump condition {\rm (\ref{RHPPSIx0b2})} for $\zeta\in\Gamma_2$}:
    We use (\ref{RHPPSIx0solution2}) to evaluate $\Psi_{\lambda,+}(\zeta)$ and
    (\ref{RHPPSIx0solution1}) to evaluate $\Psi_{\lambda,-}(\zeta)$.
    From (\ref{RHPPSIx0solution1}) and
    \cite[formulas 9.1.3, 9.1.4 and 9.6.3]{AbramowitzStegun},
    the 1,1-entry and the 2,1-entry
    on the right of (\ref{RHPPSIx0b2}) are equal to
    \begin{eqnarray}
        \nonumber
        \lefteqn{
            \Psi_{\lambda,11,-}(\zeta)+e^{-2\pi i\lambda}\Psi_{\lambda,12,-}(\zeta)
            } \\[2ex]
        \nonumber
        &=&
            \frac{1}{2}\sqrt\pi\zeta^{1/2}e^{-(\lambda+\frac{1}{4})\pi i}
            H_{\lambda+\frac{1}{2}}^{(2)}(\zeta)
            +\frac{1}{2}\sqrt\pi \zeta^{1/2}e^{-(\lambda+\frac{1}{4})\pi i}
            H_{\lambda+\frac{1}{2}}^{(1)}(\zeta) \\[2ex]
        \nonumber
        &=&
            \sqrt\pi \zeta^{1/2} e^{-(\lambda+\frac{1}{4})\pi i}
            J_{\lambda+\frac{1}{2}}(\zeta) \\[2ex]
        &=&
            \sqrt\pi \zeta^{1/2} I_{\lambda+\frac{1}{2}}(\zeta e^{-\frac{\pi i}{2}}) e^{-\frac{1}{2}\lambda\pi
            i},
    \end{eqnarray}
    and
    \begin{eqnarray}
        \nonumber
        \lefteqn{
            \Psi_{\lambda,21,-}(\zeta)+e^{-2\pi i\lambda}\Psi_{\lambda,22,-}(\zeta)
            } \\[2ex]
        \nonumber
        &=&
            \frac{1}{2}\sqrt\pi\zeta^{1/2}e^{-(\lambda+\frac{1}{4})\pi i}
            H_{\lambda-\frac{1}{2}}^{(2)}(\zeta)
            +\frac{1}{2}\sqrt\pi \zeta^{1/2}e^{-(\lambda+\frac{1}{4})\pi i}
            H_{\lambda-\frac{1}{2}}^{(1)}(\zeta) \\[2ex]
        \nonumber
        &=&
            \sqrt\pi \zeta^{1/2} e^{-(\lambda+\frac{1}{4})\pi i}
            J_{\lambda-\frac{1}{2}}(\zeta) \\[2ex]
        &=&
            -i \sqrt\pi \zeta^{1/2} I_{\lambda-\frac{1}{2}}(\zeta e^{-\frac{\pi i}{2}}) e^{-\frac{1}{2}\lambda\pi
            i},
    \end{eqnarray}
    respectively. By (\ref{RHPPSIx0solution2}) we then see that the first columns of both sides of
    (\ref{RHPPSIx0b2}) agree. From (\ref{RHPPSIx0solution1}),
    (\ref{RHPPSIx0solution2}) and
    \cite[formula 9.6.4]{AbramowitzStegun}, the second columns of both sides of (\ref{RHPPSIx0b2}) agree as well.

    \medskip

    \textit{Jump condition {\rm (\ref{RHPPSIx0b2})} for $\zeta\in\Gamma_6$}:
    We use (\ref{RHPPSIx0solution5}) to evaluate $\Psi_{\lambda,+}(\zeta)$ and
    (\ref{RHPPSIx0solution6}) to evaluate $\Psi_{\lambda,-}(\zeta)$. Since $-\zeta=\zeta e^{\pi i}$, we
    have,
    from (\ref{RHPPSIx0solution6}) and \cite[formula 9.6.4]{AbramowitzStegun},
    that the 1,2-entry and the 2,2-entry on the right of (\ref{RHPPSIx0b2}) are
    equal to
    \begin{equation}\label{Proof: PSIx0 solution: eq1}
        \Psi_{\lambda,12,-}(\zeta)
        =
            -\frac{i}{\sqrt\pi} \zeta^{1/2} e^{\frac{1}{2}\lambda\pi i} K_{\lambda+
            \frac{1}{2}}(\zeta e^{\frac{\pi i}{2}})
        =
            -\frac{i}{2} \sqrt\pi (-\zeta)^{1/2} H_{\lambda+\frac{1}{2}}^{(1)}(-\zeta)
            e^{(\lambda+\frac{1}{4})\pi i},
    \end{equation}
    and
    \begin{equation}\label{Proof: PSIx0 solution: eq2}
        \Psi_{\lambda,22,-}(\zeta)
        =
            -\frac{1}{\sqrt\pi} \zeta^{1/2} e^{\frac{1}{2}\lambda\pi i}
            K_{\lambda-\frac{1}{2}}(\zeta e^{\frac{\pi i}{2}})
        =
            \frac{i}{2} \sqrt\pi (-\zeta)^{1/2} H_{\lambda-\frac{1}{2}}^{(1)}(-\zeta)
            e^{(\lambda+\frac{1}{4})\pi i},
    \end{equation}
    respectively. So, by (\ref{RHPPSIx0solution5}) we see that the second columns of both sides
    of (\ref{RHPPSIx0b2}) agree. Since $-\zeta=\zeta e^{\pi i}$ we
    have,
    from (\ref{RHPPSIx0solution6}), (\ref{Proof: PSIx0 solution: eq1}),
    (\ref{Proof: PSIx0 solution: eq2}) and
    \cite[formulas 9.1.3, 9.1.4 and 9.6.3]{AbramowitzStegun}, that
    the 1,1-entry and the 2,1-entry on the right of (\ref{RHPPSIx0b2}) are equal to
    \begin{eqnarray}
        \nonumber
        \lefteqn{
            \Psi_{\lambda,11,-}(\zeta)+e^{-2\pi i \lambda}\Psi_{\lambda,12,-}(\zeta)
        }\\[2ex]
        \nonumber
        &=&
            - \sqrt\pi (-\zeta)^{1/2} e^{-\frac{1}{2}\lambda\pi i}
            I_{\lambda+\frac{1}{2}}(\zeta e^{\frac{\pi i}{2}})
            + \frac{1}{2} \sqrt\pi (-\zeta)^{1/2} e^{-(\lambda+\frac{1}{4})\pi i}
            H^{(1)}_{\lambda+\frac{1}{2}}(-\zeta)
             \\[2ex]
        \nonumber
        &=&
            -\sqrt\pi (-\zeta)^{1/2} e^{-(\lambda+\frac{1}{4})\pi i}
            J_{\lambda+\frac{1}{2}}(-\zeta)+ \frac{1}{2}\sqrt\pi(-\zeta)^{1/2}e^{-(\lambda+\frac{1}{4})\pi i}
            H^{(1)}_{\lambda+\frac{1}{2}}(-\zeta) \\[2ex]
        &=&
            -\frac{1}{2}\sqrt\pi (-\zeta)^{1/2} H^{(2)}_{\lambda+\frac{1}{2}}(-\zeta)e^{-(\lambda+\frac{1}{4})\pi
            i},
    \end{eqnarray}
    and
    \begin{eqnarray}
        \nonumber
        \lefteqn{
            \Psi_{\lambda,21,-}(\zeta)+e^{-2\pi i \lambda}\Psi_{\lambda,22,-}(\zeta)
        }\\[2ex]
        \nonumber
        &=&
            -i\sqrt\pi (-\zeta)^{1/2} e^{-\frac{1}{2}\lambda\pi i}
            I_{\lambda-\frac{1}{2}}(\zeta e^{\frac{\pi i}{2}})
            - \frac{1}{2} \sqrt\pi (-\zeta)^{1/2} e^{-(\lambda+\frac{1}{4})\pi i}
            H^{(1)}_{\lambda-\frac{1}{2}}(-\zeta)
             \\[2ex]
        \nonumber
        &=&
            \sqrt\pi (-\zeta)^{1/2} e^{-(\lambda+\frac{1}{4})\pi i}
            J_{\lambda-\frac{1}{2}}(-\zeta)- \frac{1}{2} \sqrt\pi (-\zeta)^{1/2} e^{-(\lambda+\frac{1}{4})\pi i}
            H^{(1)}_{\lambda-\frac{1}{2}}(-\zeta) \\[2ex]
        &=&
            \frac{1}{2}\sqrt\pi (-\zeta)^{1/2} H^{(2)}_{\lambda-\frac{1}{2}}(-\zeta)e^{-(\lambda+\frac{1}{4})\pi
            i},
    \end{eqnarray}
    respectively. By (\ref{RHPPSIx0solution5}) we then see that the first columns of both sides of
    (\ref{RHPPSIx0b2}) agree as well. We now have proven that jump condition (\ref{RHPPSIx0b2}) is
    satisfied.

    \medskip

        \textit{Jump condition {\rm (\ref{RHPPSIx0b4})}}: Similarly,
    we can prove that this jump condition is also satisfied. Here, we also use
    \cite[formula 9.1.35]{AbramowitzStegun}, and the details are left to the
    reader. This implies that the theorem is proved.
\end{proof}

Now, we explain how we get $P_{x_\nu}^{(1)}$ out of the solution
$\Psi_{\lambda_\nu}$ (depending on the parameter $\lambda_\nu$) of
the RH problem for $\Psi_{\lambda_\nu}$. We make use of the
following scalar function,
\begin{equation}\label{fx0}
    f_{x_\nu}(z)=
    \left\{\begin{array}{ll}
        i \log \varphi(z) - i \log \varphi_+(x_\nu), & \qquad \mbox{for $\Im z>0$,} \\[1ex]
        -i \log \varphi(z) - i \log \varphi_+(x_\nu), & \qquad \mbox{for $\Im z<0$,}
    \end{array}\right.
\end{equation}
which is defined and analytic for $z\in \mathbb{C}\setminus
\mathbb{R}$. For $x\in (-1,1)$ we have, since
$\varphi_+(x)\varphi_-(x)=1$, that
$f_{x_\nu,+}(x)=f_{x_\nu,-}(x)$, so that $f_{x_\nu}$ is also
analytic across the interval $(-1,1)$. The behavior of $f_{x_\nu}$
near $x_\nu$ is
\[
    f_{x_\nu}(z)=\frac{1}{\sqrt{1-x_\nu^2}}(z-x_\nu)+O\left((z-x_\nu)^2\right), \qquad \mbox{as $z\to x_\nu$.}
\]

Since $f_{x_\nu}$ is analytic near $x_\nu$, and since
$f'(x_\nu)\neq 0$, the scalar function $f_{x_\nu}$ is a one-to-one
conformal mapping on a neighborhood of $x_\nu$. So, if we choose
$\delta>0$ sufficiently small, $f_{x_\nu}$ is a one-to-one
conformal mapping on $U_{\delta,x_\nu}$ and the image of
$U_{\delta,x_\nu}$ under the mapping $\zeta=f_{x_\nu}$ is convex.

For  $x\in(-1,1)$ we have by (\ref{fx0}) that
$f_{x_\nu}(x)=\arccos(x_\nu)-\arccos x$. So, $f_{x_\nu}(x)$ is
real for $x\in(-1,1)$. If $x>x_\nu$ we have $f_{x_\nu}(x)>0$, and
if $x<x_\nu$ we have $f_{x_\nu}(x)<0$. Since $f_{x_\nu}$ is a
conformal mapping, this implies that $f_{x_\nu}$ maps
$U_{\delta,x_\nu}\cap \mathbb{C}_+$ one-to-one onto
$f_{x_\nu}(U_{\delta,x_\nu})\cap\mathbb{C}_+$, and
$U_{\delta,x_\nu}\cap \mathbb{C}_-$ one-to-one onto
$f_{x_\nu}(U_{\delta,x_\nu})\cap\mathbb{C}_-$.

\medskip

We now come back to the special choice of the contour
$\Gamma_{x_\nu}$, which we used to continuate our weight
analytically, see Section \ref{Subsection: Secont transformation}.
For $z\in\Gamma_{x_\nu}\cap\mathbb{C}_+$ we have $\arg
\varphi(z)=\arccos x_\nu$, by construction of $\Gamma_{x_\nu}$,
and for $z\in\Gamma_{x_\nu}\cap\mathbb{C}_-$ we have $\arg
\varphi(z)=-\arccos x_\nu$. By (\ref{fx0}) we then have $\Re
f_{x_\nu}(z)=0$, for $z\in\Gamma_{x_\nu}$. This implies that the
image of the contour $\Gamma_{x_\nu}$ under the mapping
$\zeta=f_{x_\nu}$ is the imaginary axis, which explains our choice
of $\Gamma_{x_\nu}$.

\begin{figure}
    \center{\resizebox{14cm}{!}{\includegraphics{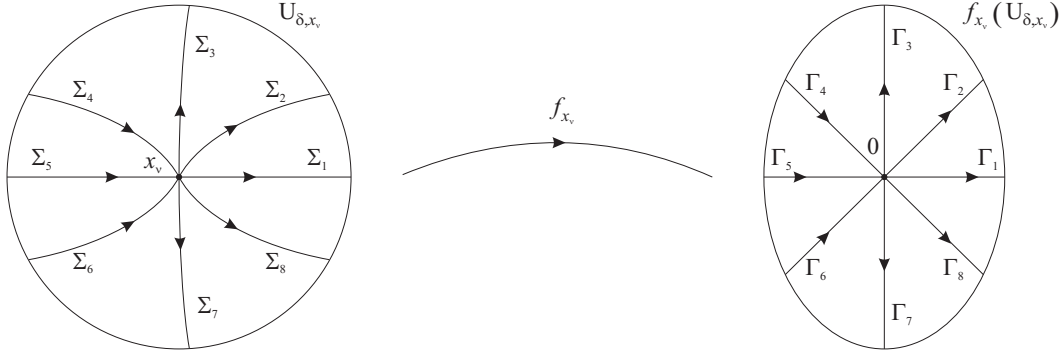}}
    \caption{The conformal mapping $f_{x_\nu}$. For every $k=1,\ldots 8$, the contour
    $\Sigma_k$ is mapped onto the part of the corresponding
    ray $\Gamma_k$ in $f_{x_\nu}(U_{\delta,x_\nu})$.}\label{figure6}}
\end{figure}

We remember that the contour $\Sigma_{x_\nu}$ was not yet
completely defined. Now, we define the contours
$\Sigma_2\cup\Sigma_4\cup\Sigma_6\cup\Sigma_8$ as the preimages of
the parts of the corresponding rays
$\Gamma_2\cup\Gamma_4\cup\Gamma_6\cup\Gamma_8$ in
$f_{x_\nu}(U_{\delta,x_\nu})$ under the mapping
$\zeta=f_{x_\nu}(z)$, see Figure \ref{figure6}. We then have
immediately that we can define
\begin{equation}
    P_{x_\nu}^{(1)}(z)=\Psi_{\lambda_\nu}(nf_{x_\nu}(z)),
\end{equation}
and $P^{(1)}_{x_\nu}$ will solve the RH problem for
$P_{x_\nu}^{(1)}$.

\begin{remark}
    We can use any one-to-one conformal mapping on
    $U_{\delta,x_\nu}$ to construct $P^{(1)}_{x_\nu}$ out of
    $\Psi_{\lambda_\nu}$. However, we have to choose it so as to
    compensate for the factor $\varphi(z)^{-n\sigma_3}$ in (\ref{P in function of
    Px01}). In the next section we will see that our choice of
    conformal mapping will do the job.
\end{remark}

\subsection{Construction of $E_{n,x_\nu}$}
    \label{Subsection: Matching condition}

We recall that for every matrix valued function $E_{n,x_\nu}$
analytic in a neighborhood of $U_{\delta,x_\nu}$, the matrix
valued function $P_{x_\nu}$ given by
\begin{equation}\label{Px0}
    P_{x_\nu}(z)=E_{n,x_\nu}(z)\Psi_{\lambda_\nu}(nf_{x_\nu}(z))W_{x_\nu}(z)^{-\sigma_3}\varphi(z)^{-n\sigma_3},
\end{equation}
satisfies conditions (a), (b) and (d) of the RH problem for
$P_{x_\nu}$. In this section we want to determine $E_{n,x_\nu}$ so
that the matching condition (c) is satisfied as well. To this end
we need to know the asymptotic behavior of $\Psi_{\lambda_\nu}$ at
infinity, and use this to determine $E_{n,x_\nu}$. At the end of
this section we also show that $E_{n,x_\nu}$ is analytic in a
neighborhood of $U_{\delta,x_\nu}$, so that the parametrix
$P_{x_\nu}$ is completely defined.

\medskip

In order to determine the asymptotic behavior of
$\Psi_{\lambda_\nu}$ at infinity, we insert the behavior of the
Bessel functions at infinity into the matrix valued function
$\Psi_{\lambda_\nu}$, given by
(\ref{RHPPSIx0solution1})--(\ref{RHPPSIx0solution8}). See
\cite[formulas 9.7.1--9.7.4]{AbramowitzStegun} for the behavior of
the modified Bessel functions at infinity, and \cite[formulas
9.2.7--9.2.10]{AbramowitzStegun} for the behavior of the Hankel
functions at infinity. Then, a straightforward calculation gives
us the asymptotic behavior of $\Psi_{\lambda_\nu}$ at infinity.
The behavior is different in each quadrant. For the upper
half-plane we find as $\zeta\to\infty$,
\begin{equation}\label{AsymptoticsPsiQuadrant1}
    \Psi_{\lambda_\nu}(\zeta)=\frac{1}{\sqrt 2}
    \begin{pmatrix}
        1 & -i \\
        -i & 1
    \end{pmatrix}
    \left[I+O\left(\frac{1}{\zeta}\right)\right]
    e^{\frac{\pi i}{4}\sigma_3} e^{-i\zeta\sigma_3} e^{-\frac{1}{2}\lambda_\nu\pi i\sigma_3},
\end{equation}
uniformly for $\zeta$ in the first quadrant, and
\begin{equation}\label{AsymptoticsPsiQuadrant2}
    \Psi_{\lambda_\nu}(\zeta)=\frac{1}{\sqrt 2}
    \begin{pmatrix}
        1 & -i \\
        -i & 1
    \end{pmatrix}
    \left[I+O\left(\frac{1}{\zeta}\right)\right]
    e^{\frac{\pi i}{4}\sigma_3} e^{-i\zeta\sigma_3} e^{\frac{1}{2}\lambda_\nu\pi i\sigma_3},
\end{equation}
uniformly for $\zeta$ in the second quadrant. For the lower
half-plane we find as $\zeta\to\infty$,
\begin{equation}\label{AsymptoticsPsiQuadrant3}
    \Psi_{\lambda_\nu}(\zeta)=\frac{1}{\sqrt 2}
    \begin{pmatrix}
        1 & -i \\
        -i & 1
    \end{pmatrix}
    \left[I+O\left(\frac{1}{\zeta}\right)\right]
    e^{\frac{\pi i}{4}\sigma_3} e^{-i\zeta\sigma_3} e^{\frac{1}{2}\lambda_\nu\pi i \sigma_3}
    \begin{pmatrix}
        0 & -1 \\
        1 & 0
    \end{pmatrix},
\end{equation}
uniformly for $\zeta$ in the third quadrant, and
\begin{equation}\label{AsymptoticsPsiQuadrant4}
    \Psi_{\lambda_\nu}(\zeta)=\frac{1}{\sqrt 2}
    \begin{pmatrix}
        1 & -i \\
        -i & 1
    \end{pmatrix}
    \left[I+O\left(\frac{1}{\zeta}\right)\right]
    e^{\frac{\pi i}{4}\sigma_3} e^{-i\zeta\sigma_3} e^{-\frac{1}{2}\lambda_\nu\pi i \sigma_3}
    \begin{pmatrix}
        0 & -1 \\
        1 & 0
    \end{pmatrix},
\end{equation}
uniformly for $\zeta$ in the fourth quadrant.

\medskip

Now, we use the asymptotic behavior
(\ref{AsymptoticsPsiQuadrant1})--(\ref{AsymptoticsPsiQuadrant4})
of $\Psi_{\lambda_\nu}$ at infinity to determine $E_{n,x_\nu}$. We
explain this only for the region $\partial U_{\delta,x_\nu}\cap
K_{x_\nu}^{\rm I}$. The other cases are similar and the details
are left to the reader. For $z\in
\partial U_{\delta,x_\nu}\cap K_{x_\nu}^{\rm I}$ we have, since $f_{x_\nu}$ is a one-to-one
conformal mapping on $U_{\delta,x_\nu}$, that $nf_{x_\nu}(z)$ lies
in the first quadrant, cf.\! Figure \ref{figure6}. So, we may use
(\ref{AsymptoticsPsiQuadrant1}) to evaluate the asymptotic
behavior of $\Psi_{\lambda_\nu}(nf_{x_\nu}(z))$ as $n\to\infty$.
Since $\Im z>0$, we have by (\ref{fx0}),
\[
    e^{-inf_{x_\nu}(z)}=\varphi_+(x_\nu)^{-n}\varphi(z)^{n}.
\]
Using (\ref{Px0}) and (\ref{AsymptoticsPsiQuadrant1}) we then find
\begin{eqnarray}
    \nonumber
    P_{x_\nu}(z)N^{-1}(z)
    &=&
        E_{n,x_\nu}(z)\frac{1}{\sqrt 2}
        \begin{pmatrix}
            1 & -i \\
            -i & 1
        \end{pmatrix}
        \left[I+O\left(\frac{1}{n}\right)\right] \\[1ex]
    \nonumber
    &&
        \qquad\qquad\times\,
        e^{\frac{\pi i}{4}\sigma_3}
        \varphi_+(x_\nu)^{-n\sigma_3}e^{-\frac{1}{2}\lambda_\nu\pi
        i \sigma_3}W_{x_\nu}(z)^{-\sigma_3}N^{-1}(z),
\end{eqnarray}
as $n\to\infty$, uniformly for $z\in \partial U_{\delta,x_\nu}\cap
K_{x_\nu}^{\rm I}$. So, in order that the matching condition is
satisfied we define for $z\in U\cap K_{x_\nu}^{\rm I}$,
\begin{equation}
    E_{n,x_\nu}(z)=N(z) W_{x_\nu}(z)^{\sigma_3} e^{\frac{1}{2}\lambda_\nu\pi i \sigma_3} \varphi_+(x_\nu)^{n\sigma_3}
    e^{-\frac{\pi i}{4}\sigma_3} \frac{1}{\sqrt 2}
    \begin{pmatrix}
        1 & i \\
        i & 1
    \end{pmatrix}.
\end{equation}
\begin{remark}
    With this $E_{n,x_\nu}$ we see that
    \begin{eqnarray*}
        P_{x_\nu}(z)N^{-1}(z)&=&
            N(z)W_{x_\nu}(z)^{\sigma_3}\varphi_+(x_\nu)^{n\sigma_3}
            \left[I+O\left(\frac{1}{n}\right)\right] \\[1ex]
        &&
            \qquad\qquad
            \times\,
            \varphi_+(x_\nu)^{-n\sigma_3}W_{x_\nu}(z)^{-\sigma_3}N^{-1}(z),
    \end{eqnarray*}
    as $n\to\infty$, uniformly for $z\in\partial U_{\delta,x_\nu}\cap K_{x_\nu}^{\rm I}$.
    Since $|\varphi_+(x_\nu)|=1$,
    and since $W_{x_\nu}$ as well as all entries of $N$ are
    bounded and bounded away from 0 on $\partial U_{\delta,x_\nu}$,
    the matching condition is satisfied on $\partial U_{\delta,x_\nu}\cap K_{x_\nu}^{\rm I}$.
\end{remark}
Similarly, we use (\ref{fx0}) and (\ref{Px0}) together with the
asymptotic behavior
(\ref{AsymptoticsPsiQuadrant2})--(\ref{AsymptoticsPsiQuadrant4})
of $\Psi_{\lambda_\nu}$ at infinity to determine $E_{n,x_\nu}$ in
the other regions. Straightforward calculations then show that we
have to define $E_{n,x_\nu}(z)$ for $z\in
U\setminus(\mathbb{R}\cup\Gamma_{x_\nu})$ as
\begin{equation}\label{Enx0}
    E_{n,x_\nu}(z)= E_\nu(z) \varphi_+(x_\nu)^{n\sigma_3} e^{-\frac{\pi i}{4}\sigma_3} \frac{1}{\sqrt 2}
    \begin{pmatrix}
        1 & i \\
        i & 1
    \end{pmatrix},
\end{equation}
where the matrix valued function $E_\nu(z)$ does not depend on
$n$, is analytic for $z\in
U\setminus(\mathbb{R}\cup\Gamma_{x_\nu})$ and given by
\begin{eqnarray}
    \label{NTildeQuadrant1}
    E_\nu(z)
    &=&
        N(z) W_{x_\nu}(z)^{\sigma_3} e^{\frac{1}{2}\lambda_\nu\pi
        i\sigma_3}, \qquad \mbox{for $z\in U\cap K_{x_\nu}^{\rm I}$,} \\[2ex]
    \label{NTildeQuadrant2}
    E_\nu(z)
    &=&
        N(z) W_{x_\nu}(z)^{\sigma_3} e^{-\frac{1}{2}\lambda_\nu\pi
        i\sigma_3}, \qquad \mbox{for $z\in U\cap K_{x_\nu}^{\rm II}$,} \\[2ex]
    \label{NTildeQuadrant3}
    E_\nu(z)
    &=&
        N(z) W_{x_\nu}(z)^{\sigma_3}
        \begin{pmatrix}
            0 & 1 \\
            -1 & 0
        \end{pmatrix}
        e^{-\frac{1}{2}\lambda_\nu\pi i\sigma_3}, \qquad \mbox{for $z\in U\cap K_{x_\nu}^{\rm III}$,} \\[2ex]
    \label{NTildeQuadrant4}
    E_\nu(z)
    &=&
        N(z) W_{x_\nu}(z)^{\sigma_3}
        \begin{pmatrix}
            0 & 1 \\
            -1 & 0
        \end{pmatrix}
        e^{\frac{1}{2}\lambda_\nu\pi i\sigma_3}, \qquad \mbox{for $z\in U\cap K_{x_\nu}^{\rm IV}$.}
\end{eqnarray}

\medskip

Now, everything is fine, except for the fact that $E_{n,x_\nu}$ is
analytic in $U\setminus(\mathbb{R}\cup \Gamma_{x_\nu})$, but we
want it to be analytic in a full neighborhood of $x_\nu$. This
will be proven in the next proposition.

\begin{proposition}\label{Proposition E constant}
    The matrix valued function $E_{n,x_\nu}$ defined by {\rm (\ref{Enx0})}--{\rm
    (\ref{NTildeQuadrant4})} is analytic in $U\setminus( (-\infty,x_{\nu-1}]\cup[x_{\nu+1},\infty))$.
\end{proposition}

\begin{proof}
    By (\ref{Enx0}) it suffices to prove that $E_\nu$ is
    analytic in $U\setminus( (-\infty,x_{\nu-1}]\cup[x_{\nu+1},\infty))$. We will check that
    $E_\nu$ has no jumps on
    $(x_{\nu-1},x_{\nu+1})\setminus\{x_\nu\}$ and
    $(\Gamma_{x_\nu}\cap U)\setminus\{x_\nu\}$, and in addition
    that the isolated singularity of $E_\nu$ at $x_\nu$ is removable. Let
    $(x_{\nu-1},x_{\nu+1})$ be oriented from the left to the right,
    and let $\Gamma_{x_\nu}\cap U$ be oriented so that it points away from $x_\nu$,
    cf.\! Figure \ref{figure4}.

    \medskip

    For $x\in (x_{\nu-1},x_\nu)$ we use (\ref{NTildeQuadrant2}) to
    evaluate $E_{\nu,+}(x)$ and (\ref{NTildeQuadrant3}) to
    evaluate $E_{\nu,-}(x)$. From (\ref{Wx0^2}) we have
    $W_{x_\nu,+}(x)=w(x)^{1/2}e^{\lambda_\nu\pi i}$ and
    $W_{x_\nu,-}(x)=w(x)^{1/2}e^{-\lambda_\nu\pi i}$. Therefore, by
    (\ref{RHPNb}), (\ref{NTildeQuadrant2}) and
    (\ref{NTildeQuadrant3}),
    \begin{eqnarray}
        \nonumber
        E_{\nu,+}(x)
        &=&
            N_-(x)
            \begin{pmatrix}
                0 & w(x) \\
                -w(x)^{-1} & 0
            \end{pmatrix}
            w(x)^{\sigma_3/2} e^{\lambda_\nu\pi i\sigma_3} e^{-\frac{1}{2}\lambda_\nu\pi i\sigma_3} \\[1ex]
        \nonumber
        &=&
            N_-(x)w(x)^{\sigma_3/2}e^{-\lambda_\nu\pi i\sigma_3}
            \begin{pmatrix}
                0 & 1 \\
                -1 & 0
            \end{pmatrix}
            e^{-\frac{1}{2}\lambda_\nu\pi i\sigma_3} \\[1ex]
        \nonumber
        &=&
            E_{\nu,-}(x).
    \end{eqnarray}
    Hence, $E_{\nu}$ is analytic across $(x_{\nu-1},x_\nu)$.
    Similarly, we have by (\ref{RHPNb}), (\ref{Wx0^2}), (\ref{NTildeQuadrant1}) and
    (\ref{NTildeQuadrant4}) that
    $E_{\nu}$ is analytic across $(x_\nu,x_{\nu+1})$ as well.

    For $z\in (\Gamma_{x_\nu}\cap U)\cap\mathbb{C}_+$ we
    use (\ref{NTildeQuadrant2}) to evaluate $E_{\nu,+}(z)$ and
    (\ref{NTildeQuadrant1}) to evaluate $E_{\nu,-}(z)$. From
    (\ref{Wx0+Wx0-invers}) we have $W_{x_\nu,+}(z)W_{x_\nu,-}(z)^{-1}=e^{\lambda_\nu\pi
    i}$. Therefore, by
    (\ref{NTildeQuadrant1}) and (\ref{NTildeQuadrant2}),
    \begin{eqnarray}
        \nonumber
        E_{\nu,+}(z)
        &=&
            N(z)W_{x_\nu,+}(z)^{\sigma_3}
            e^{-\frac{1}{2}\lambda_\nu\pi i\sigma_3} \\[1ex]
        \nonumber
        &=&
            E_{\nu,-}(z)
            e^{-\frac{1}{2}\lambda_\nu\pi i\sigma_3}
            W_{x_\nu,-}(z)^{-\sigma_3}W_{x_\nu,+}(z)^{\sigma_3}
            e^{-\frac{1}{2}\lambda_\nu\pi i\sigma_3} \\[1ex]
        \nonumber
        &=&
            E_{\nu,-}(z),
    \end{eqnarray}
    so that $E_\nu$ is analytic across
    $(\Gamma_{x_\nu}\cap U)\cap\mathbb{C}_+$. Similarly, we have by (\ref{Wx0+Wx0-invers}), (\ref{NTildeQuadrant3})
    and (\ref{NTildeQuadrant4}) that
    $E_\nu$ is analytic across $(\Gamma_{x_\nu}\cap U)\cap\mathbb{C}_-$ as well. We thus have proven that
    $E_\nu$ is analytic in $U\setminus((-\infty,x_{\nu-1}]\cup[x_{\nu+1},\infty)\cup\{x_\nu\})$.

    \medskip

    It remains to prove that the isolated singularity of $E_\nu$ at $x_\nu$ is
    removable. We have by
    (\ref{Szegofunction}) and (\ref{Wx0}) that $D(z)$ behaves like $c_1
    |z-x_\nu|^{\lambda_\nu}$
    and $W_{x_\nu}(z)$ like $c_2
    |z-x_\nu|^{\lambda_\nu}$ as $z\to x_\nu$, where $c_1$
    and $c_2$ are non-zero constants. Therefore
    \[
        \frac{W_{x_\nu}(z)}{D(z)}=O(1),\qquad \mbox{and}\qquad \frac{D(z)}{W_{x_\nu}(z)}=O(1),
        \qquad\mbox{as $z\to x_\nu$.}
    \]
    So, by (\ref{RHPNsolution}), each entry of
    $N(z)W_{x_\nu}(z)^{\sigma_3}$ remains bounded as $z\to x_\nu$.
    This implies by (\ref{NTildeQuadrant1})--(\ref{NTildeQuadrant4}) that each entry
    of $E_\nu$ remains bounded as $z\to x_\nu$, so that the
    isolated singularity of $E_\nu$ at $x_\nu$ is removable. Therefore, the
    proposition is proved.
\end{proof}
This ends the construction of the local parametrix near $x_\nu$.

\medskip

We recall that we also wanted the local parametrix $P_{x_\nu}$ to
be invertible, see Remark \ref{Remark: invers exists}. We will
show that
\begin{equation}\label{det Px0=1}
    \det P_{x_\nu}\equiv 1.
\end{equation}
This is analogous as in \cite[Section 7]{KMVV} and we will just
give a sketch of the proof. Since $E_{n,x_\nu}$ is a product of
four matrices all with determinant 1, see
(\ref{Enx0})--(\ref{NTildeQuadrant4}), it suffices to prove from
(\ref{Px0}) that $\det \Psi_{\lambda_\nu}\equiv 1$. Using part (b)
of the RH problem for $\Psi_{\lambda_\nu}$ we find that $\det
\Psi_{\lambda_\nu}$ is analytic in $\mathbb{C}\setminus\{0\}$. If
we then use the behavior of $\Psi_{\lambda_\nu}$ near $0$ stated
in part (c) of the RH problem the isolated singularity of $\det
\Psi_{\lambda_\nu}$ at $0$ has to be removable, so that $\det
\Psi_{\lambda_\nu}$ is an entire function. Using the asymptotics
of $\Psi_{\lambda_\nu}$ at infinity given by
(\ref{AsymptoticsPsiQuadrant1})--(\ref{AsymptoticsPsiQuadrant4})
we have that $\det \Psi_{\lambda_\nu}(\zeta)\to 1$ as
$\zeta\to\infty$. By Liouville's theorem we then have that $\det
\Psi_{\lambda_\nu}\equiv 1$, so that also $\det P_{x_\nu}\equiv
1$.

\section{Asymptotics of the recurrence coefficients}

In this section we will determine a complete asymptotic expansion
of the recurrence coefficients $a_n$ and $b_n$ as $n\to\infty$.
Recall that $a_n$ and $b_n$ have been formulated in terms of the
solution of the RH problem for $Y$, see (\ref{an in Y}) and
(\ref{bn in Y}). The asymptotic analysis of the RH problem for $Y$
has been done in Section 3, and unfolding the series of
transformations $Y\mapsto T\mapsto S\mapsto R$, see \cite[Section
9]{KMVV} for details, we find
\begin{equation}\label{an in function of R}
    a_n^2=\lim_{z\to\infty}\left(-\frac{D_\infty^2}{2i}+zR_{12}(z;n,w)\right)
        \left(zR_{21}(z;n,w)+\frac{1}{2iD_\infty^2}\right),
\end{equation}
and
\begin{equation}\label{bn in function of R}
    b_n=\lim_{z\to\infty}\left(z-zR_{11}(z;n+1,w)R_{22}(z;n,w)\right).
\end{equation}

\begin{remark}
    We note that, see \cite[Lemma 8.3]{KMVV},
    \begin{equation}\label{result Lemma 8.3 KMVV}
        \left\| R(z)-I \right\| =
        O\left(\frac{1}{n|z|}\right),\qquad \mbox{as $n\to\infty$,}
    \end{equation}
    uniformly for $|z|\geq 2$, where $\|\cdot\|$ is any matrix
    norm. Inserting this into (\ref{an in function of R}) and (\ref{bn in function of R})
    we find the known asymptotic behavior of the recurrence
    coefficients, cf. \cite{Golinskii},
    \[
        a_n=\frac{1}{2}+O(1/n),\qquad
        b_n=O(1/n),
        \qquad \mbox{as $n\to\infty$.}
    \]
    In the rest of the paper we will develop the $O(1/n)$ terms
    into complete asymptotic expansions in powers of $1/n$.
\end{remark}

In order to determine a complete asymptotic expansion of $a_n$ and
$b_n$ we will work as follows. In Section \ref{Subsection:
Asymptotic expansion for Delta}, we will determine a complete
asymptotic expansion of the jump matrix for $R$ in powers of $1/n$
as $n\to\infty$. As a result, we obtain in Section
\ref{Subsection: Asymptotic expansion for R} a complete asymptotic
expansion of $R$. The coefficients in this expansion can be
calculated explicitly via residue calculus, and we will determine
the order $1/n$ term. Finally, in Section \ref{Subsection: Proof
of our result} we will use this to prove Theorem \ref{Theorem:
asymptotic expansion recurrence coefficients}.

\subsection{Asymptotic expansion of $\Delta$}
    \label{Subsection: Asymptotic expansion for Delta}

Denote the jump matrix for $R$ as $I+\Delta$. Then, from condition
(c) of the RH problem for $R$,
\begin{eqnarray}
    \label{Delta(s)1}
    \Delta(z)
    &=&
        P_{x_\nu}(z)N^{-1}(z)-I, \qquad \mbox{for $z\in \partial
            U_{\delta,x_\nu}$ and $\nu=0,\ldots ,n_0+1$,} \\[2ex]
    \label{Delta(s)2}
    \Delta(z)
    &=&
        N(z)
        \begin{pmatrix}
            1 & 0\\
            w(z)^{-1}\varphi(z)^{-2n} & 1
        \end{pmatrix}
        N^{-1}(z)-I, \qquad \mbox{for $z\in \Sigma_R\setminus
            (\cup_{\nu=0}^{n_0+1}\partial U_{\delta,x_\nu})$.}
\end{eqnarray}
In this section we will show that $\Delta$ has an asymptotic
expansion in powers of $1/n$ of the form
\begin{equation}\label{asymptotic expansion Delta}
    \Delta(z)\sim \sum_{k=1}^\infty
    \frac{\Delta_k(z,n)}{n^k},\qquad \mbox{as $n\to\infty$,}
\end{equation}
uniformly for $z\in\Sigma_R$, and we will also determine the
coefficients $\Delta_k(z,n)$ explicitly.

\begin{remark}
    The $n$-dependance of the coefficients in the expansion
    will come from the factor $\varphi_+(x_\nu)^{n\sigma_3}$ in the parametrices
    near the algebraic singularities $x_\nu$.
\end{remark}

On the lips of the lens, $\Delta$ vanishes at an exponential rate,
cf.\! \cite[Section 7]{KMVV}. This implies for every $k$,
\begin{equation}
    \Delta_k(z,n)= 0,\qquad \mbox{for $z\in \Sigma_R\setminus
            (\cup_{\nu=0}^{n_0+1}\partial U_{\delta,x_\nu})$.}
\end{equation}
On the circles near $\pm 1$, the asymptotic expansion
(\ref{asymptotic expansion Delta}) of $\Delta$ is known, see
\cite[Section 8]{KMVV}. The restriction of $\Delta_k(\cdot,n)$ to
$\partial U_{\delta,-1}\cup
\partial U_{\delta,1}$ is given by \cite[(8.5) and (8.6)]{KMVV}
and does not depend on $n$. It has a meromorphic continuation to
$U_{\delta_0,-1}\cup U_{\delta_0,1}$ for some $\delta_0>\delta$,
with poles of order at most $[(k+1)/2]$ at $\pm 1$. For details we
refer to \cite[Section 8]{KMVV}.

\medskip

So, it remains to determine the asymptotic expansion of $\Delta$
on the circles near the algebraic singularities. Fix
$\nu\in\{1,\ldots ,n_0\}$. By (\ref{Px0}), (\ref{Enx0}) and
(\ref{Delta(s)1}), we have for $z\in\partial U_{\delta,x_\nu}$,
\begin{eqnarray}
    \nonumber
    \Delta(z)&=&
        E_\nu(z)\varphi_+(x_\nu)^{n\sigma_3}e^{-\frac{\pi
        i}{4}\sigma_3}\frac{1}{\sqrt 2}
        \begin{pmatrix}
            1 & i \\
            i & 1
        \end{pmatrix} \\[1ex]
    \label{Asymptotic expansion Delta near x0 eq1}
    && \qquad
        \times\,
        \Psi_{\lambda_\nu}(nf_{x_\nu}(z))\varphi(z)^{-n\sigma_3}W_{x_\nu}(z)^{-\sigma_3}N^{-1}(z)-I.
\end{eqnarray}
Here, the matrix valued function $\Psi_{\lambda_\nu}$ is
constructed out of Bessel functions, which have a complete
asymptotic expansion at infinity. This implies that
$\Psi_{\lambda_\nu}(nf_{x_\nu}(z))$ also has a complete asymptotic
expansion as $n\to\infty$. Inserting the asymptotic expansions of
the modified Bessel functions at infinity \cite[formulas
9.7.1--9.7.4]{AbramowitzStegun} into (\ref{RHPPSIx0solution1}),
and the asymptotic expansions of the Hankel functions
\cite[formulas 9.2.7--9.2.10]{AbramowitzStegun} into
(\ref{RHPPSIx0solution2}), we obtain
\begin{eqnarray}
    \nonumber
    \Psi_{\lambda_\nu}(\zeta)&\sim&
        \frac{1}{\sqrt 2}
            \begin{pmatrix}
                1 & -i \\
                -i & 1
            \end{pmatrix}
            \left[
            I+\sum_{k=1}^\infty \frac{i^k}{2^{k+1} \zeta^k}
            \begin{pmatrix}
                (-1)^k s_{\lambda_\nu,k} & -i t_{\lambda_\nu,k} \\[1ex]
                i (-1)^k t_{\lambda_\nu,k} & s_{\lambda_\nu,k}
            \end{pmatrix}\right] \\[2ex]
        \label{Asymptotic expansion Delta near x0 eq2}
        && \qquad \times\,
            e^{\frac{\pi i}{4}\sigma_3}e^{-\frac{1}{2}\lambda_\nu\pi i\sigma_3}e^{-i\zeta\sigma_3},
\end{eqnarray}
as $\zeta\to\infty$, uniformly for $\zeta$ in the first quadrant.
Here, the constants $s_{\lambda_\nu,k}$ and $t_{\lambda_\nu,k}$
are given by
\begin{equation}\label{sgammatgamma}
        s_{\lambda_\nu,k}=(\lambda_\nu+\frac{1}{2},k)+(\lambda_\nu-\frac{1}{2},k),\qquad
        t_{\lambda_\nu,k}=(\lambda_\nu+\frac{1}{2},k)-(\lambda_\nu-\frac{1}{2},k),
\end{equation}
where
\[
    (\nu,k)=\frac{(4\nu^2-1)(4\nu^2-9)\ldots(4\nu^2-(2k-1)^2)}{2^{2k}k!}.
\]
For $z\in \partial U_{\delta,x_\nu}\cap K_{x_\nu}^{\rm I}$  we
have, since $f_{x_\nu}$ is a one-to-one conformal mapping on
$U_{\delta,x_\nu}$ that $nf_{x_\nu}(z)$ lies in the first
quadrant, see Figure \ref{figure6}. So, we may use
(\ref{Asymptotic expansion Delta near x0 eq2}) to determine the
asymptotic expansion of $\Psi_{\lambda_\nu}(nf_{x_\nu}(z))$ as
$n\to\infty$. Since $\Im z>0$ we have by (\ref{fx0}),
\[
    e^{-inf_{x_\nu}(z)}=\varphi_+(x_\nu)^{-n}\varphi(z)^{n}.
\]
Therefore, by (\ref{Asymptotic expansion Delta near x0 eq2}),
\begin{eqnarray*}
    \lefteqn{
    \frac{1}{\sqrt 2}
    \begin{pmatrix}
        1 & i\\
        i & 1
    \end{pmatrix}
    \Psi_{\lambda_\nu}(nf_{x_\nu}(z))\varphi(z)^{-n\sigma_3}\sim}
    \\[2ex]
    && \qquad \left[
            I+\sum_{k=1}^\infty \frac{i^k}{2^{k+1} f_{x_\nu}(z)^k n^k}
            \begin{pmatrix}
                (-1)^k s_{\lambda_\nu,k} & -i t_{\lambda_\nu,k} \\[1ex]
                i (-1)^k t_{\lambda_\nu,k} & s_{\lambda_\nu,k}
            \end{pmatrix}\right]e^{\frac{\pi i}{4}\sigma_3}\varphi_+(x_\nu)^{-n\sigma_3}
            e^{-\frac{1}{2}\lambda_\nu\pi i\sigma_3},
\end{eqnarray*}
as $n\to\infty$, uniformly for $z\in \partial U_{\delta,x_\nu}\cap
K_{x_\nu}^{\rm I}$. Inserting this into (\ref{Asymptotic expansion
Delta near x0 eq1}) we have by (\ref{NTildeQuadrant1}), and the
fact that $\varphi_+(x_\nu)^n$ remains bounded and bounded away
from 0 as $n\to\infty$ (since $|\varphi_+(x_\nu)|=1$),
\begin{eqnarray}
    \nonumber
    \Delta(z) &\sim &
        \sum_{k=1}^\infty \frac{i^k}{2^{k+1} f_{x_\nu}(z)^k}
        E_\nu(z)\varphi_+(x_\nu)^{n\sigma_3}
        \begin{pmatrix}
            (-1)^k s_{\lambda_\nu,k} & - t_{\lambda_\nu,k} \\[1ex]
            - (-1)^k t_{\lambda_\nu,k} & s_{\lambda_\nu,k}
        \end{pmatrix} \\[2ex]
    && \qquad\times\,
        \varphi_+(x_\nu)^{-n\sigma_3} E_\nu^{-1}(z)\frac{1}{n^k},
\end{eqnarray}
as $n\to\infty$, uniformly for $z\in \partial U_{\delta,x_\nu}\cap
K_{x_\nu}^{\rm I}$. Similarly, we find the same asymptotic
expansion on the other regions of $\partial U_{\delta,x_\nu}$. The
details are left to the reader. Thus, for $z\in\partial
U_{\delta,x_\nu}$ the coefficients of the expansion
(\ref{asymptotic expansion Delta}) for $\Delta$ are given by
\begin{equation}\label{DeltakUdeltax0}
        \Delta_k(z,n)=
            \frac{i^k}{2^{k+1} f_{x_\nu}(z)^k}
            E_\nu(z)\varphi_+(x_\nu)^{n\sigma_3}
            \begin{pmatrix}
                (-1)^k s_{\lambda_\nu,k} & - t_{\lambda_\nu,k} \\[1ex]
                - (-1)^k t_{\lambda_\nu,k} & s_{\lambda_\nu,k}
            \end{pmatrix}
            \varphi_+(x_\nu)^{-n\sigma_3}E_\nu^{-1}(z).
\end{equation}
\begin{remark}
    These coefficients depend on $n$ through the factors $\varphi_+(x_\nu)^{\pm
    n\sigma_3}$. Since $|\varphi_+(x_\nu)|=1$, the coefficients $\Delta_k(z,n)$ remain bounded as $n\to\infty$,
    which is necessary to get an
    asymptotic expansion of the form (\ref{asymptotic expansion Delta}).
\end{remark}

We note that $f_{x_\nu}^{-k}$ is analytic in
$\mathbb{C}\setminus((-\infty,-1]\cup[1,\infty))$ except for a
pole of order $k$ at $x_\nu$, see the discussion at the end of
Section \ref{Subsection: Construction Px01}. From the proof of
Proposition \ref{Proposition E constant} and the fact that $\det
E_\nu=1$ we have that $E_\nu$ as well as $E_\nu^{-1}$ are analytic
in $U\setminus((-\infty,x_{\nu-1}]\cup[x_{\nu+1},\infty))$. So,
the restriction of $\Delta_k(\cdot,n)$ to $\partial
U_{\delta,x_\nu}$ has a meromorphic continuation to a neighborhood
$U_{\delta_0,x_\nu}$ of $x_\nu$ for some $\delta_0>\delta$, with a
pole of order $k$ at $x_\nu$.

\subsection{Asymptotic expansion of $R$}
    \label{Subsection: Asymptotic expansion for R}

We recall that $\Delta$ possesses an asymptotic expansion in
powers of $1/n$ of the form (\ref{asymptotic expansion Delta})
with oscillatory terms in the expansion. Following the argument
that leads to \cite[(4.115)]{DKMVZ1}, this implies that $R$ itself
possesses an asymptotic expansion in powers of $1/n$ given by
\begin{equation}\label{asymptotic expansion R}
    R(z;n,w)\sim I+\sum_{k=1}^\infty\frac{R_k(z,n)}{n^k},\qquad\mbox{as $n\to\infty$,}
\end{equation}
uniformly for $z\in\mathbb{C}\setminus\Sigma_R$. Here, for every
$k$ and $n$,
\begin{equation}\label{additive RH problem eq1}
    \mbox{$R_k(\cdot,n)$ is analytic in $\mathbb{C}\setminus(\cup_{\nu=0}^{n_0+1}\partial
    U_{\delta,x_\nu})$,}
\end{equation}
and
\begin{equation}\label{additive RH problem eq2}
    R_k(z,n)=O(1/z),\qquad \mbox{as $z\to\infty$.}
\end{equation}
The $n$-dependance in the coefficients $R_k(z,n)$ arises through
the oscillatory terms in the expansion of $\Delta$.

\medskip

We will now determine, similar as in \cite[Section 8]{KMVV}, the
coefficient $R_1(z,n)$ explicitly. Expanding the jump relation
$R_+=R_-(I+\Delta)$, and collecting the terms with $1/n$ we have
\begin{equation}\label{additive RH problem eq3}
    R_{1,+}(s,n)-R_{1,-}(s,n)=\Delta_1(s,n),\qquad
    \mbox{for $s\in\cup_{\nu=0}^{n_0+1}\partial
    U_{\delta,x_\nu}$,}
\end{equation}
which is, together with (\ref{additive RH problem eq1}) and
(\ref{additive RH problem eq2}) an additive RH problem. This can
easily be solved using the Sokhotskii-Plemelj formula, but in our
case we can write down an explicit solution as follows. Since
$\Delta_1(z,n)$ is analytic in neighborhoods of $z=\pm 1$ and
$z=x_\nu$ for $\nu=1,\ldots n_0$, except for simple poles at those
points, see Section \ref{Subsection: Asymptotic expansion for
Delta}, we can write
\[
    \Delta_1(z,n)=\frac{A^{(1)}(n)}{z-1}+O(1),\quad\mbox{as $z\to
    1$,}\qquad
    \Delta_1(z,n)=\frac{B^{(1)}(n)}{z+1}+O(1),\quad\mbox{as $z\to
    -1$,}
\]
and
\[
    \Delta_1(z,n)=\frac{C_\nu^{(1)}(n)}{z-x_\nu}+O(1),\quad\mbox{as $z\to
    x_\nu$,}
\]
for certain constant matrices $A^{(1)}(n), B^{(1)}(n)$ and
$C_\nu^{(1)}(n)$.
\begin{remark}
    Since $\Delta_1(s,n)$ is independent of $n$ for $s\in\partial U_{\delta,-1}\cup\partial
    U_{\delta,1}$, see Section \ref{Subsection: Asymptotic expansion for
    Delta},
    the residues $A^{(1)}(n)$ and $B^{(1)}(n)$ of $\Delta_1(z,n)$ at $z=1$ and $z=-1$, respectively, are
    also independent of $n$. The $n$-dependance of the residue
    $C_\nu^{(1)}(n)$ at $z=x_\nu$
    follows from the oscillatory terms $\varphi_+(x_\nu)^{\pm n\sigma_3}$ in
    $\Delta_1(s,n)$ near $x_\nu$, see (\ref{DeltakUdeltax0}).
\end{remark}
By inspection we then see that
\begin{equation}
    R_1(z,n)=\left\{\begin{array}{l}
        {\ds
        \frac{A^{(1)}(n)}{z-1}+\frac{B^{(1)}(n)}{z+1}+\sum_{\nu=1}^{n_0}\frac{C_\nu^{(1)}(n)}{z-x_\nu}},
        \qquad \mbox{for $z\in\mathbb{C}\setminus(
            \cup_{\nu=0}^{n_0+1}
            U_{\delta,x_\nu})$,}\\[5ex]
        {\ds
        \frac{A^{(1)}(n)}{z-1}+\frac{B^{(1)}(n)}{z+1}+\sum_{\nu=1}^{n_0}\frac{C_\nu^{(1)}(n)}{z-x_\nu}-\Delta_1(z,n)},
        \qquad \mbox{for $z\in \cup_{\nu=0}^{n_0+1}
            U_{\delta,x_\nu}$,}
    \end{array}\right.
\end{equation}
satisfies the additive RH problem (\ref{additive RH problem
eq1})--(\ref{additive RH problem eq3}). So, we need to determine
the constant matrices $A^{(1)}(n), B^{(1)}(n)$ and
$C_\nu^{(1)}(n)$ for $\nu=1,\ldots ,n_0$. For $A^{(1)}(n)$ and
$B^{(1)}(n)$ we have found, see \cite[Section 8]{KMVV},
\begin{equation}\label{A(1)}
    A^{(1)}(n)=\frac{4\alpha^2-1}{16}D_\infty^{\sigma_3}
    \begin{pmatrix}
        -1 & i\\
        i & 1
    \end{pmatrix}
    D_\infty^{-\sigma_3},
\end{equation}
\begin{equation}\label{B(1)}
    B^{(1)}(n)=\frac{4\beta^2-1}{16}D_\infty^{\sigma_3}
    \begin{pmatrix}
        1 & i\\
        i & -1
    \end{pmatrix}
    D_\infty^{-\sigma_3},
\end{equation}
which is clearly independent of $n$. It remains to determine the
residue $C_\nu^{(1)}(n)$ of $\Delta_1(z,n)$ at $z=x_\nu$, for
every $\nu=1,\ldots ,n_0$. Fix $\nu\in\{1,\ldots ,n_0\}$. Since
$E_\nu(z)$ and $E_\nu^{-1}(z)$ are analytic in a neighborhood of
$z=x_\nu$, and since
\[
    f_{x_\nu}(z)^{-1}=\sqrt{1-x_\nu^2}\frac{1}{z-x_\nu}+O(1), \qquad
    \mbox{as $z\to x_\nu$,}
\]
we have by (\ref{sgammatgamma}) and (\ref{DeltakUdeltax0}) that
the residue $C_\nu^{(1)}(n)$ of $\Delta_1(z,n)$ at $z=x_\nu$ is
given by
\begin{equation}\label{determination C1n eq1}
    C_\nu^{(1)}(n)=\frac{i}{4}\sqrt{1-x_\nu^2}
    E_\nu(x_\nu)\varphi_+(x_\nu)^{n\sigma_3}
    \begin{pmatrix}
        -2\lambda_\nu^2 & -2\lambda_\nu \\
        2\lambda_\nu & 2\lambda_\nu^2
    \end{pmatrix}
    \varphi_+(x_\nu)^{-n\sigma_3}E_\nu^{-1}(x_\nu).
\end{equation}

We want to simplify this expression. So, we need to find
convenient expressions for $E_\nu(x_\nu)$ and $
E_\nu^{-1}(x_\nu)$, and substitute these into (\ref{determination
C1n eq1}). Since $E_\nu$ is analytic near $x_\nu$ we determine
$E_\nu(x_\nu)$ by the following limit
\[
    E_\nu(x_\nu)=\lim_{x\downarrow x_\nu;\, x\in\mathbb{R}}E_{\nu}(x)=
    \lim_{x\downarrow x_\nu;\, x\in\mathbb{R}}E_{\nu,+}(x).
\]
Here, we take the limit from $x$ to $x_\nu$ on the real axis from
the right. The last equality follows from the fact that $E_\nu$
has no jumps on $(x_\nu,x_{\nu+1})$, see Proposition
\ref{Proposition E constant}. From (\ref{RHPNsolution}) and
(\ref{NTildeQuadrant1}) we then find
\begin{equation}\label{Ntilde x0 eq1}
    E_\nu(x_\nu)=\lim_{x\downarrow x_\nu;\, x\in\mathbb{R}}D_\infty^{\sigma_3}
        \begin{pmatrix}
            \frac{a_+(x)+a_+(x)^{-1}}{2} &
                \frac{a_+(x)-a_+(x)^{-1}}{2i}\\[1ex]
            \frac{a_+(x)-a_+(x)^{-1}}{-2i} & \frac{a_+(x)+a_+(x)^{-1}}{2}
        \end{pmatrix}
        \left(\frac{W_{x_\nu,+}(x)}{D_+(x)}\right)^{\sigma_3}e^{\frac{1}{2}\lambda_\nu\pi
        i\sigma_3}.
\end{equation}
The Szeg\H{o} function satisfies
$D_+(x)=\sqrt{w(x)}e^{-i\psi_\nu(x)}$ for $x\in(x_\nu,x_{\nu+1})$,
see Lemma \ref{Lemma: D+}, where $\psi_\nu$ is given by
(\ref{psi_nu}). By (\ref{Wx0^2}) we have
$W_{x_\nu,+}(x)=\sqrt{w(x)}e^{-\lambda_\nu\pi i}$, so that by
(\ref{Theorem: Phi}) and (\ref{psi_nu})
\[
    \lim_{x\downarrow x_\nu;\, x\in\mathbb{R}}\frac{W_{x_\nu,+}(x)}{D_+(x)}e^{\frac{1}{2}\lambda_\nu\pi i}
        =e^{i\psi_\nu(x_\nu)}e^{-\frac{1}{2}\lambda_\nu\pi
        i}=e^{-\frac{1}{2}\Phi_\nu i}.
\]
Inserting this into (\ref{Ntilde x0 eq1}) and using the following
identities, which hold for $x\in(-1,1)$,
\[
    \frac{a_+(x)+a_+(x)^{-1}}{2}=\frac{e^{-\frac{\pi i}{4}}}{\sqrt
    2(1-x^2)^{1/4}}\varphi_+(x)^{1/2},
\]
\[
    \frac{a_+(x)-a_+(x)^{-1}}{2i}=\frac{e^{-\frac{\pi i}{4}}}{\sqrt
    2(1-x^2)^{1/4}}i\varphi_+(x)^{-1/2},
\]
we then find
\begin{equation}\label{Ntilde x0}
    E_\nu(x_\nu)=\frac{e^{-\frac{\pi i}{4}}}{\sqrt
    2(1-x_\nu^2)^{1/4}}D_\infty^{\sigma_3}
    \begin{pmatrix}
        \varphi_+(x_\nu)^{1/2} & i\varphi_+(x_\nu)^{-1/2}\\
        -i\varphi_+(x_\nu)^{-1/2} & \varphi_+(x_\nu)^{1/2}
    \end{pmatrix}
    e^{-\frac{1}{2}\Phi_\nu i\sigma_3}.
\end{equation}
Taking inverse, we find
\begin{equation}\label{NtildeInver x0}
    E_\nu^{-1}(x_\nu)=\frac{e^{-\frac{\pi i}{4}}}{\sqrt
    2(1-x_\nu^2)^{1/4}}e^{\frac{1}{2}\Phi_\nu i\sigma_3}
    \begin{pmatrix}
        \varphi_+(x_\nu)^{1/2} & -i\varphi_+(x_\nu)^{-1/2}\\
        i\varphi_+(x_\nu)^{-1/2} & \varphi_+(x_\nu)^{1/2}
    \end{pmatrix}
    D_\infty^{-\sigma_3}.
\end{equation}

Now, we insert (\ref{Ntilde x0}) and (\ref{NtildeInver x0}) into
(\ref{determination C1n eq1}). Using the identity
$\varphi_+(x_\nu)=\exp(i\arccos x_\nu)$ we then find after a
straightforward calculation that the residue $C^{(1)}_\nu(n)$ of
$\Delta_1(z,n)$ at $z=x_\nu$ is given by
\begin{equation}\label{Cnu1(n)}
    C_\nu^{(1)}(n)=D_\infty^{\sigma_3}
        \begin{pmatrix}
            C_{\nu,11}(n) & C_{\nu,12}(n)\\
            C_{\nu,21}(n) & -C_{\nu,11}(n)
        \end{pmatrix}D_\infty^{-\sigma_3},
\end{equation}
where
\begin{eqnarray}\label{Cnu1(n)11}
    C_{\nu,11}(n)&=& -\frac{1}{2}\lambda_\nu^2 x_\nu+\frac{1}{2}\lambda_\nu\sin\bigl(2n\arccos
    x_\nu-\Phi_\nu\bigr),
\end{eqnarray}
\begin{eqnarray}
    \nonumber
    C_{\nu,12}(n) &=&
        \frac{i}{2}\lambda_\nu^2-\frac{i}{2}\lambda_\nu
        x_\nu\sin\bigl(2n\arccos x_\nu-\Phi_\nu\bigr)\\[1ex]
    \label{Cnu1(n)12}
    && \qquad\qquad -\,
        \frac{i}{2}\lambda_\nu\sqrt{1-x_\nu^2}\cos\bigl(2n\arccos
        x_\nu-\Phi_\nu\bigr),
\end{eqnarray}
\begin{eqnarray}
    \nonumber
    C_{\nu,21}(n) &=&
        \frac{i}{2}\lambda_\nu^2-\frac{i}{2}\lambda_\nu
        x_\nu\sin\bigl(2n\arccos x_\nu-\Phi_\nu\bigr)\\[1ex]
    \label{Cnu1(n)21}
    && \qquad\qquad +\,
        \frac{i}{2}\lambda_\nu\sqrt{1-x_\nu^2}\cos\bigl(2n\arccos
        x_\nu-\Phi_\nu\bigr).
\end{eqnarray}
This ends the determination of $R_1(z,n)$.

For general $k$, we get that $R_k(z,n)$ in the region
$\mathbb{C}\setminus(\cup_{\nu=0}^{n_0+1}U_{\delta,x_\nu})$ is a
rational function with poles at $\pm 1$ and at the algebraic
singularities $x_\nu$. The residues at 1 and $-1$ are denoted by
$A^{(k)}(n)$ and $B^{(k)}(n)$ respectively and may depend on $n$.
The residue at every $x_\nu$ depends on $n$ and is denoted by
$C^{(k)}_\nu(n)$. We then get
\[
    R_k(z,n)=\frac{A^{(k)}(n)}{z-1}+\frac{B^{(k)}(n)}{z+1}+\sum_{\nu=1}^{n_0}\frac{C^{(k)}_\nu(n)}{z-x_\nu}+
    O(1/z^2), \qquad
    \mbox{as $z\to\infty$.}
\]
The residues $A^{(k)}(n), B^{(k)}(n)$ and $C^{(k)}_\nu(n)$ can be
determined in a similar fashion, but for our purpose it suffices
to know $R_1(z,n)$.

\subsection{Proof of Theorem \ref{Theorem: asymptotic expansion recurrence coefficients}}
    \label{Subsection: Proof of our result}

We are now ready to determine a complete asymptotic expansion of
the recurrence coefficients $a_n$ and $b_n$. The idea is to insert
the asymptotic expansion (\ref{asymptotic expansion R}) of $R$
into (\ref{an in function of R}) and (\ref{bn in function of R}).

\begin{varproof}{\bf of Theorem \ref{Theorem: asymptotic expansion recurrence
coefficients}.} We recall that, see (\ref{an in function of R}),
\[
    a_n^2=\lim_{z\to\infty}\left(-\frac{D_\infty^2}{2i}+zR_{12}(z;n,w)\right)
        \left(zR_{21}(z;n,w)+\frac{1}{2iD_\infty^2}\right).
\]
We may take the limit $z\to\infty$ in the asymptotic expansion
(\ref{asymptotic expansion R}) of $R$, cf.\! \cite[Section
9]{KMVV}, to obtain
\begin{eqnarray}
    \nonumber
    a_n^2 &\sim&
        \left(-\frac{D_\infty^2}{2i}+
        \sum_{k=1}^\infty\frac{A^{(k)}_{12}(n)+B^{(k)}_{12}(n)+\sum_{\nu=1}^{n_0}C_{\nu,12}^{(k)}(n)}{n^k}\right)
        \\[2ex]
    \label{Expansionansquare}
    && \qquad\qquad\times\,
        \left(\sum_{k=1}^\infty\frac{A^{(k)}_{21}(n)+B^{(k)}_{21}(n)
        +\sum_{\nu=1}^{n_0}C_{\nu,21}^{(k)}(n)}{n^k}+\frac{1}{2i
        D_\infty^2}\right),
\end{eqnarray}
as $n\to\infty$. Expanding this we find a complete asymptotic
expansion of $a_n^2$, and this leads to a complete asymptotic
expansion of $a_n$ in powers of $1/n$ as $n\to\infty$. By
(\ref{A(1)}), (\ref{B(1)}), (\ref{Cnu1(n)}), (\ref{Cnu1(n)12}),
(\ref{Cnu1(n)21}) and (\ref{Expansionansquare}) the first terms in
the asymptotic expansion of $a_n^2$ are
\begin{eqnarray}
    \nonumber
    a_n^2 &=&
        \frac{1}{4}+\frac{1}{2i}\left[D_\infty^{-2}\left(A_{12}^{(1)}(n)+B_{12}^{(1)}(n)+
        \sum_{\nu=1}^{n_0}C_{\nu,12}^{(1)}(n)\right)\right. \\[1ex]
    \nonumber
    && \qquad\qquad\qquad -\, \left.
        D_\infty^2\left(A_{21}^{(1)}(n)+B_{21}^{(1)}(n)+\sum_{\nu=1}^{n_0}C_{\nu,21}^{(1)}(n)\right)\right]\frac{1}{n}
        +O\left(\frac{1}{n^2}\right)
        \\[2ex]
    \label{Expansionansquarefirstterms}
    &=&
        \frac{1}{4}-\sum_{\nu=1}^{n_0}\frac{\lambda_\nu}{2}\sqrt{1-x_\nu^2}\cos\bigl(2n\arccos
        x_\nu-\Phi_\nu\bigr)\frac{1}{n}+O\left(\frac{1}{n^2}\right),
\end{eqnarray}
as $n\to\infty$. From this we then get, after a simple
calculation, that the coefficient with the $1/n$ term in the
asymptotic expansion of $a_n$ is given by (\ref{Theorem: A1(n)}).
So, the statements about the recurrence coefficient $a_n$ are
proved.

\medskip
Similarly, we can prove the statements about the recurrence
coefficient $b_n$. If we take in (\ref{bn in function of R}) the
limit $z\to\infty$ in the asymptotic expansion (\ref{asymptotic
expansion R}) of $R$, cf.\! \cite[Section 9]{KMVV}, we find
\begin{eqnarray}
    \nonumber
    b_n &\sim&
        \lim_{z\to\infty} -z\left(\sum_{k=1}^\infty \frac{(R_k)_{11}(z,n+1)}{(n+1)^k}
        +\sum_{k=1}^\infty \frac{(R_k)_{22}(z,n)}{n^k}\right)\\[2ex]
    \nonumber
    &=&
        -\sum_{k=1}^\infty\left(\frac{A^{(k)}_{11}(n+1)+B^{(k)}_{11}(n+1)
        +\sum_{\nu=1}^{n_0}C_{\nu,11}^{(k)}(n+1)}{(n+1)^k}\right. \\[2ex]
    \label{Expansionbn}
    &&
        \left.
        \qquad\qquad\qquad
        +\,\frac{A^{(k)}_{22}(n)+B^{(k)}_{22}(n)+\sum_{\nu=1}^{n_0}C_{\nu,22}^{(k)}(n)}{n^k}\right),
\end{eqnarray}
as $n\to\infty$. From this we get a complete asymptotic expansion
of $b_n$ in powers of $1/n$, and by (\ref{A(1)}), (\ref{B(1)}),
(\ref{Cnu1(n)}), (\ref{Cnu1(n)11}) and (\ref{Expansionbn}) the
coefficient with the $1/n$ term in the asymptotic expansion of
$b_n$ is given by
\begin{eqnarray}
    \nonumber
    B_1(n) &=& -\left(A^{(1)}_{11}(n+1)+A^{(1)}_{22}(n)\right)
        -\left(B^{(1)}_{11}(n+1)+B^{(1)}_{22}(n)\right) \\[1ex]
    \nonumber
    &&
        \qquad\qquad\qquad
        -\, \sum_{\nu=1}^{n_0}\left(C_{\nu,11}^{(1)}(n+1)+C_{\nu,22}^{(1)}(n)\right) \\[1ex]
    \nonumber
    &=&
        -\sum_{\nu=1}^{n_0}\frac{\lambda_\nu}{2}\left[\sin\Bigl((2n+1)\arccos
        x_\nu-\Phi_\nu+\arccos x_\nu\Bigr)\right. \\[1ex]
    \nonumber
    &&
        \qquad\qquad\qquad\left. -\,
        \sin\Bigl((2n+1)\arccos x_\nu-\Phi_\nu-\arccos x_\nu\Bigr)\right]
        \\[1ex]
    &=&
        - \sum_{\nu=1}^{n_0}\lambda_\nu\sqrt{1-x_\nu^2}\cos\bigl((2n+1)\arccos
        x_\nu-\Phi_\nu\bigr).
\end{eqnarray}
Therefore, the theorem is proved.
\end{varproof}

\subsection*{Acknowledgements}
    I thank my advisor Arno Kuijlaars for useful discussions, good ideas, and his great
    support.

\end{document}